\documentclass[twoside]{article}
\usepackage{latexsym}
\usepackage{amsfonts}
\usepackage{graphicx}

\newtheorem{theorem}{Theorem}[section]
\newtheorem{lemma}{Lemma}[section]

\newtheorem{definition}{Definition}[section]
\newtheorem{prop}{Proposition}[section]
\newcommand{\qed}{\hfill$\Box$\par\medskip}

\newenvironment{Proof}{\noindent{\sc Proof.}}{\qed}

\def\bhag#1{\noindent
\setcounter{equation}{0}
\section{#1}
}
\def\bfgk#1{{{#1}\kern-5.5pt{#1}}}

\def\HH{{\mathbf H}}

\def\RR{{\mathbb R}}

\def\PPI{{{\rm I}\kern-1pt\Pi}}
\def\SS{{\mathbb S}}

\def\b #1;{{\bf #1}}
\def\x{{\bf x}}

\def\y{{\bf y}}

\def\e{\epsilon}

\def\O{{\cal O}}

\def\C{{\mathcal C}}

\def\W{{\bf W}}
\def\WW{{\mathbb W}}

\def\ip#1#2{{\langle {#1}, {#2}\rangle}}

\def\esssup{\mathop{\hbox{{\rm ess sup}}}}

\def\be{\begin{equation}}
\def\ee{\end{equation}}
\def\bea{\begin{eqnarray}}
\def\eea{\end{eqnarray}}
\def\eref#1{(\ref{#1})}
\def\disp{\displaystyle}

\def\prob{\mbox{{\rm Prob }}}
\def\donchitre#1#2{\vskip 6.5cm\noindent
\parbox[t]{1in}{\special{eps:#1.eps x=6.5cm y=5.5cm}}
\hbox to 7cm{}\parbox[t]{0.0cm}{\special{eps:#2.eps x=6.5cm y=5.5cm}}}

\title{Localized linear polynomial operators and  quadrature formulas on the sphere}
\author{Q. T. Le Gia\thanks{The research of this author was supported by Australian Research Council under its Centres of
Excellence Program.}\\
School of Mathematics,  University of New South Wales, \\
Sydney, NSW 2052,    Australia \\
qlegia@maths.unsw.edu.au\\ [3ex]
H.~N.~Mhaskar\thanks{The research of this author was supported, in part, by grant DMS-0605209 from the National Science Foundation and  grant W911NF-04-1-0339 from the U.S. Army Research Office.}\\
Department of Mathematics, California
State University\\
Los Angeles, California, 90032, U.S.A.\\
hmhaska@calstatela.edu}

\date{}

\graphicspath{{converted_graphics/}}
\begin{document}
\maketitle
\hspace{1.0cm}

\begin{abstract}
The purpose of this paper is to construct universal, auto--adaptive, localized, linear, polynomial  (-valued) operators based on scattered data on the (hyper--)sphere $\SS^q$ ($q\ge 2$). The approximation and localization  properties of our operators are studied  theoretically in deterministic as well as probabilistic settings. Numerical experiments are presented to demonstrate their superiority over traditional least squares and discrete Fourier projection polynomial approximations. An essential ingredient in our construction is the construction of quadrature formulas based on scattered data, exact for integrating spherical polynomials of (moderately) high degree.  Our formulas are based on scattered sites; i.e., in contrast to such well known formulas as Driscoll--Healy formulas, we need not choose the location of the sites in any particular manner. While the previous attempts to construct such formulas have yielded formulas exact for spherical polynomials of degree at most $18$, we are able to construct formulas exact for spherical polynomials of degree $178$.  
\end{abstract}

\medskip\noindent
\textbf{Keywords:} Quadrature formulas, localized kernels, polynomial quasi--interpolation, learning theory on the sphere.\\
\textbf{AMS classification:} 65D32, 41A10, 41A25\\

\bhag{Introduction}
The problem of approximation  of functions on the sphere arises in almost all applications involving modeling of data collected on the surface of the earth. More recent applications such as manifold matching and neural networks lead to the approximation of functions on the unit sphere $\SS^q$ embedded in the Euclidean space $\RR^{q+1}$ for integers $q\ge 3$ as well. Various applications in learning theory, meteorology, cosmology, and geophysics require analysis of 
\emph{scattered data} collected on the sphere \cite{freedensurvey, freedenbook1, freedenbook2}. This means that the data is of the form $\{(\xi, f(\xi))\}$ for some unknown function $f :\SS^q\to\RR$, where one has no control on the choice of the sites $\xi$. 

There are many methods to model such data: spherical splines, radial basis functions (called zonal function networks in this context), etc. However, the most traditional method is to approximate by spherical polynomials; i.e., restrictions of  algebraic polynomials in $q+1$ variables to $\SS^q$. Apart from tradition, some important advantages of polynomials are that they are eigenfunctions of many pseudo--differential operators which arise in practical applications, and that they are infinitely smooth. Unlike in the case of spline approximation with a given degree of the piecewise component polynomials, global polynomial approximation does not exhibit a saturation property \cite[Section~2, Chapter~11]{devlorbk}; i.e., for an arbitrary sequence $\delta_n\downarrow 0$, it is possible to find a continuous function on the sphere, not itself a polynomial, which can be approximated by spherical polynomials of degree at most $n$ uniformly within $\delta_n$, $n\ge 1$.   In \cite{mnw2, zfquadpap}, we have shown how a good polynomial approximation yields also a  good zonal function network approximation. In \cite{sbfconvpap}, we have shown that the approximation spaces determined by zonal function network approximation are the same as those determined by polynomial approximations.

\begin{figure}[h]
\centering
\begin{minipage}{0.4\textwidth}
  \centering
  \includegraphics[width=\textwidth,keepaspectratio]{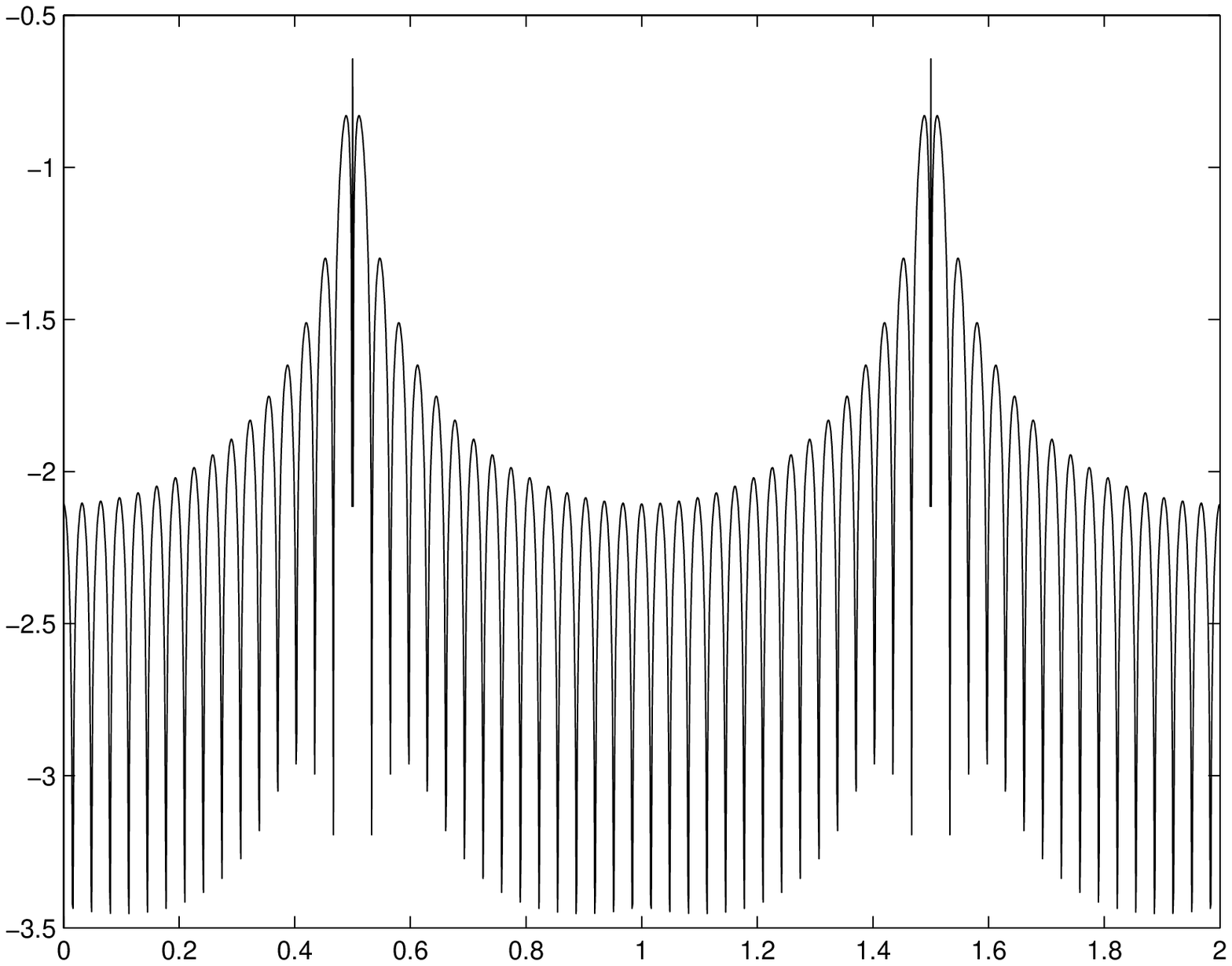}
  \end{minipage}
\begin{minipage}{0.4\textwidth}

  \centering
  \includegraphics[width=\textwidth,keepaspectratio]{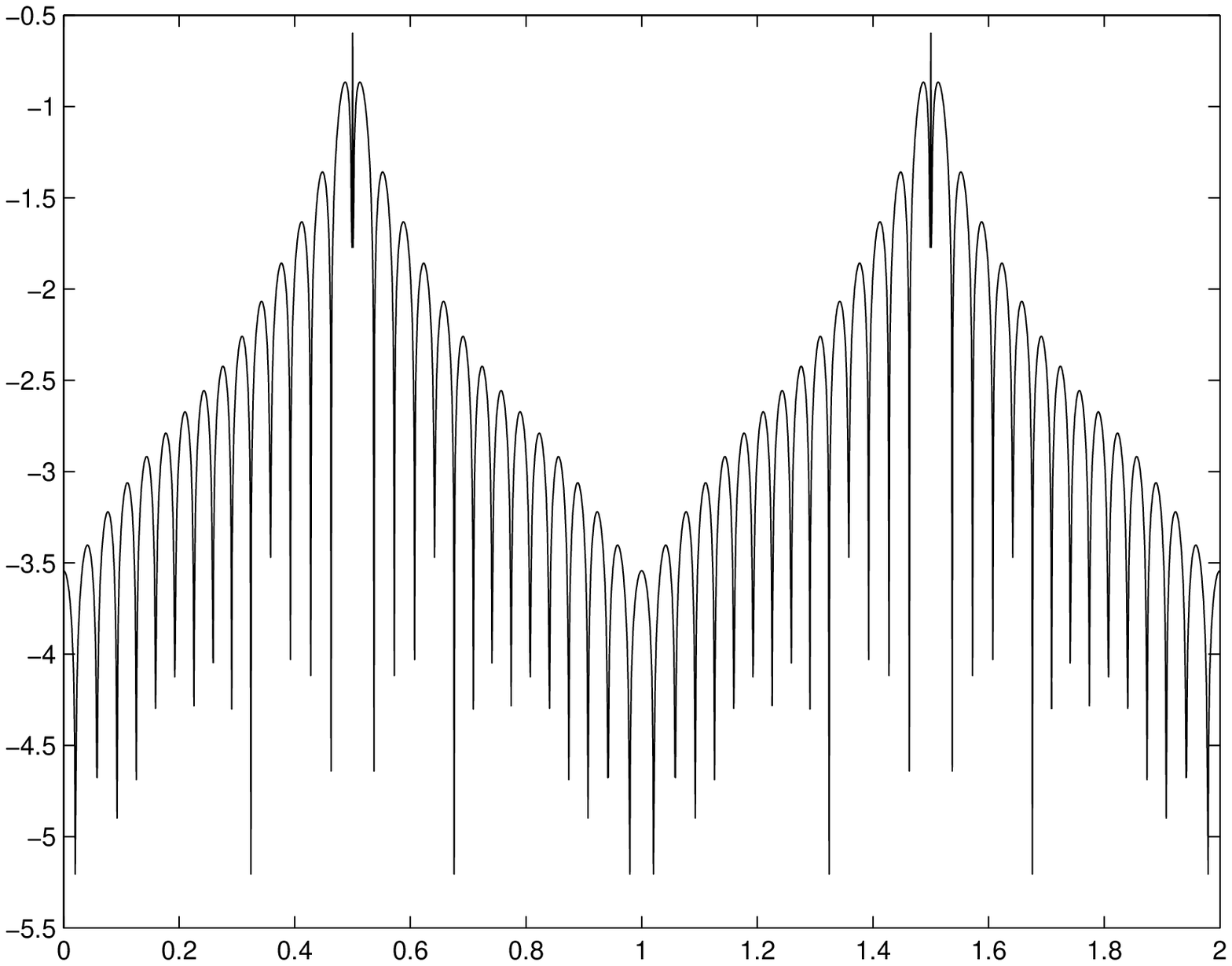}
  \end{minipage}

\caption{The log--plot of the absolute error between the function $x\mapsto |\cos x|^{1/4}$, and (left) its Fourier projection (right) trigonometric polynomial obtained by our summability operator, where the Fourier coefficients are estimated by 128 point FFT. The order of the trigonometric polynomials is $31$ in each case. The numbers on the $x$ axis are in multiples of $\pi$, the actual absolute errors are $10^{y}$.}
\label{trigfig}
\end{figure}

To illustrate the issues to be discussed in this paper,  we consider an example in the case $q=1$, or equivalently, the case of $2\pi$--periodic functions on the real line. In this discussion only, let $f(x)=|\cos x|^{1/4}$, $x\in\RR$. In Figure~\ref{trigfig}(left), we show the log-plot of the absolute errors between $f$ and its (trigonometric) Fourier projection of order $31$, where the Fourier coefficients are estimated by a $128$ point DFT. In Figure~\ref{trigfig}(right), we show a similar log--plot where the Fourier projection is replaced by a suitable summability operator (described more precisely in \eref{sigmastardef}), yielding again a trigonometric polynomial of order $31$. It is clear that our summability operator is far more localized than the Fourier projection; i.e., the error in approximation decreases more rapidly as one goes away from the singularities at $\pi/2$ and $3\pi/2$. The maximum error on $[3\pi/4,5\pi/4]$ is $0.0103$ for the projection,  $0.0028$ for our operator. Out of the 2048 points considered for the test,  the error by the summability operator is less than $10^{-3}$ at $38.96\%$ points, the corresponding percentage for the projection is only $4.88\%$. In contrast to free--knot spline approximation, our summability operator is universal; i.e., its construction (convolution with a kernel) does not require any a priori knowledge about the location of singularities of the target function. It yields a single, globally defined trigonometric polynomial, computed using global data. Nevertheless, it is auto--adaptive, in the sense that the error in approximation on different subintervals adjusts itself according to the smoothness of the target function on these subintervals. In \cite{trigwave, loctrigwav}, we have given a very detailed analysis of the approximation properties of these operators in the case $q=1$.  

Our computation based on a $128$ point DFT implies that the values of the function are available at 128 equidistant points. If only a scattered data is available, the following method is often used (especially in the context of approximation on the sphere) to estimate the values needed for the DFT. For each point $\xi$, we consider the nearest point of the form $2\pi k/128$, and imagine that the value of $f$ at this point is $f(\xi)$, taking averages in the case of multiplicities, and interpolating in the case of gaps. If we use our summability operator, estimating the Fourier coefficients in this way, then the  maximum error on $[3\pi/4,5\pi/4]$ is $0.0357$, and the proportion of points where the error is less than $10^{-3}$ is  $7.08\%$.  It is clear that a careful construction of quadrature formulas is essential to obtain good approximation results.

The purpose of this paper is to construct universal, auto--adaptive, localized, linear, polynomial  (-valued) operators based on scattered data on  $\SS^q$ ($q\ge 2$) and to analyse their approximation properties.  An essential ingredient in our construction is the construction of quadrature formulas based on scattered data, exact for integrating spherical polynomials of (moderately) high degree, and satisfying certain technical conditions known as the Marcinkiewicz--Zygmund (M--Z) conditions. Our construction is different from the usual construction of quadrature formulas (designs) studied in numerical analysis, where one has a choice of the placement of nodes. In \cite{mnw1}, we had proved the existence of such quadrature formulas for scattered data. However, previous efforts to compute such formulas did not yield exactness beyond degree $18$ polynomials. This was a severe limitation on the practical applications of our theoretical constructions. We will show that a very simple idea of solving a system of equations involving a Gram matrix yields surprisingly good results, in particular, quadrature formulas exact for integrating polynomials of degree as high as $178$. Gram matrices are typically ill--conditioned. However, we will show both theoretically and numerically that the ones which we use are, in fact, very well conditioned.  We will introduce another algorithm of theoretical interest to compute data dependent orthogonal polynomials, and use these to compute the quadrature formulas in a memory efficient manner. To the best of our knowledge, this is the first effort to extend the univariate constructions in Gautschi's book \cite{gautschi} to a multivariate setting. Considering that computation of classical spherical harmonics is a very delicate task, requiring many tricks based on the special function properties of these polynomials for a stable computation, it is not expected that our computation of data dependent orthogonal polynomials with no such special function properties would be stable. In describing this algorithm, we hope to stimulate further research in this interesting direction. We note that even if this algorithm is not as stable for high degrees as the other algorithm, it yields satisfactory quadrature formulas exact for integrating polynomials of degree $32$. Most importantly, our new found ability to compute quadrature formulas for moderately high degrees allows us to offer our operators as a viable, practical method of approximation, even superior to the commonly used methods of least squares and Fourier projection as far as localized approximation is concerned. 

An additional problem is when the available values of the target function are noisy. One may assume that the noise is an additive random variable with mean zero. It is also routine in learning theory to assume that the random variables have a bounded range. This assumption is usually satisfied with a high probability even if the random variables do not actually have a bounded range. However, one does not typically know the actual distribution of these random variables. We obtain probabilistic estimates in this setting on the global and local approximations by our operators. To underline the practical utility of our operators, we use them for modelling the MAGSAT data supplied to us by Dr. Thorsten Maier, obtaining results comparable to those obtained by other techniques.

In Section~\ref{backsect}, we review certain facts about spherical polynomials, the existence of quadrature formulas to integrate these, a few properties of the quadrature weights,  and certain polynomial kernels which we will need throughout the paper. In Section~\ref{operatorsect}, we study the approximation properties of the linear polynomial operators.  The new results here are Theorems~\ref{probopertheoa} and \ref{localapproxtheo}. The first parts of these theorems were proved essentially in \cite{locsmooth}, but not stated in the form given here. In order to apply these operators in practice, one needs quadrature formulas exact for high degree spherical polynomials. Explicit algorithms to construct such formulas are described in Section~\ref{quadconsect}.   The new results in this section  are Theorems~\ref{lsqwttheo} and \ref{recurtheo}. Numerical results are presented in Section~\ref{numsect}, and the proofs of all new results  are given in Section~\ref{proofsect}. The paper is a result of a long process, involving discussions with a number of mathematicians. In particular, it is our pleasure to acknowledge the support and encouragement of Mahadevan Ganesh,  Thorsten Maier, Volker Michel, Dominik Michel, Ian Sloan, and Joe Ward. We are also grateful to the two referees and Fred Hickernell for their many useful suggestions for the improvement of the first draft of this paper.
  
\bhag{Background}\label{backsect}
In this section, we review some known results regarding spherical polynomials and localized polynomial kernels.
\subsection{Spherical polynomials}\label{sphpolybacksect}
Let $q\ge 1$ be an integer, $\SS^q$ be the unit sphere  embedded in the Euclidean space $\RR^{q+1}$; i.e.,\\
$
\SS^q:=\{(x_1,\ldots,x_{q+1})\in \RR^{q+1}\ :\ \sum_{k=1}^{q+1} x_k^2=1\},
$
and $\mu_q$ be its Lebesgue surface measure, normalized so that $\mu_q(\SS^q)=1$. The surface area of $\SS^q$ is $\disp \frac{2\pi^{(q+1)/2}}{\Gamma((q+1)/2)}$.
For $\delta>0$, a spherical cap with radius $\delta$ and center $\x_0\in\SS^q$ is defined by
$$
\SS^q_\delta(\x_0):=\{ \x\in\SS^q\ :\ \arccos(\x\cdot\x_0) \le \delta\}.
$$

If $1\le p\le \infty$,  and $f: \SS^q\to\RR$ is  measurable, we write
$$
\|f\|_{p} :=\cases{\{\int_{\SS^q} |f(\x)|^pd\mu_q(\x)\}^{1/p}, & if $1\le p <\infty$,\cr
                                       \esssup_{\x\in \SS^q} |f(\x)|, & if $p=\infty$.\cr} 
$$
The space of all Lebesgue measurable functions on $\SS^q$ such that $\|f\|_{p} <\infty$ will be
denoted by $L^p$, with the usual convention that two functions are considered equal as
elements of this space if they are equal almost everywhere.
The symbol $C(\SS^q)$ denotes the class of all  continuous,  real valued functions on $\SS^q$, equipped with the norm $\|\circ\|_{\infty}$. 

For a real number $x\ge 0$, let $\Pi_x^q$ denote the class of all spherical polynomials of degree at most $x$. (This is the same as the class $\Pi_n^q$, where $n$ is the largest integer not exceeding $x$. However, our extension of the notation allows us, for example, to use the simpler notation $\Pi_{n/2}^q$ rather than the more cumbersome notation $\Pi_{\lfloor n/2\rfloor}^q$.)
For a fixed integer $\ell\ge 0$, the restriction to $\SS^q$ of a
homogeneous harmonic polynomial of exact degree $\ell$ is called a spherical
harmonic of degree $\ell$. Most of the following information is based
on \cite{Muller}, \cite[Section~IV.2]{SteinWeiss}, and \cite[Chapter XI]{bateman}, although we use a
different notation. The class of all spherical harmonics of degree
$\ell$ will be denoted by $\HH^q_\ell$. The
spaces $\HH^q_\ell$ are mutually orthogonal relative to
the inner product of $L^2$. For any integer $n\ge 0$, we have $\Pi^q_n = \bigoplus_{\ell=0
}^{n}\HH^q_\ell$. The dimension of $\HH^q_\ell$ is given by
\be\label{hkdim}
d\,^q_\ell := \dim \HH^q_\ell= \left\{
\begin{array}{cl}
\disp{\frac{2\ell+q-1}{\ell+q-1} {\ell+q-1 \choose \ell}}, & \mbox{if
}\ell\ge 1,\\[3ex]
1,& \mbox{if }\ell=0.
\end{array}
\right.
\ee
and that of $\Pi^q_n$ is $\sum_{\ell=0}^n d\,^q_\ell=d\,^{q+1}_n$. Furthermore,
$L^2=\mbox{\rm $L^2$--closure}\big\{\bigoplus_{\ell=0}^\infty \HH^q_\ell\big\}$.
Hence, if we choose an orthonormal basis
$\{Y_{\ell,k}\,:\,k=1,\ldots, d^q_\ell\}$ for each $\HH^q_\ell$, then
the set $\{Y_{\ell,k}\,:\,\ell=0,1,\ldots\, \mbox{and }k=1,\ldots,
d\,^q_\ell\}$ is a complete orthonormal basis for $L^2$. One has the
well-known addition formula \cite{Muller} and \cite[Chapter XI, Theorem 4]{bateman}: 
\begin{equation}
\label{addformula}
 \sum_{k=1}^{d\,^q_\ell} Y_{\ell,k}(\b x;){Y_{\ell,k}(\zeta)} =
\frac{2^{q-1}\Gamma(q/2)^2}{\Gamma(q)}p_\ell(1)p_\ell(\x\cdot\zeta), \qquad
\ell=0,1,\cdots, 
\end{equation}
where $p_\ell:=p_\ell^{(q/2-1,q/2-1)}$ is the orthonormalized Jacobi polynomial with positive leading coefficient: 
$$
\int_{-1}^1 p_\ell(t)p_k(t)(1-t^2)^{q/2-1} dt =\left\{\begin{array}{ll}
1, & \mbox{ if $\ell=k$, }\\
0, &\mbox{otherwise}.
\end{array}\right.
$$
In particular, for $\x\in\SS^q$, $\ell=0,1,\cdots$,
\be\label{constchristoffel}
\sum_{k=1}^{d_\ell^q} Y_{\ell,k}^2(\x)=\frac{2^{q-1}\Gamma(q/2)^2}{\Gamma(q)}p_\ell(1)^2=\int_{\SS^q}\sum_{k=1}^{d_\ell^q} Y_{\ell,k}^2(\zeta)d\mu_q(\zeta)=d_\ell^q.
\ee

\subsection{Localized polynomial kernels}
Let $h: [0,\infty)\to\RR$ be a compactly supported function, and $t>0$. We define for $u\in\RR$,
\be\label{kerndef}
\Phi_t(h;u):=\frac{2^{q-1}\Gamma(q/2)^2}{\Gamma(q)}\sum_{\ell=0}^\infty h(\ell/t) p_\ell(1)p_\ell(u)
\ee
and define $\Phi_t(h;u)=0$ if $t\le 0$. 

 In the sequel, we adopt the following convention
regarding constants. The letters $c, c_1,\cdots$ will denote generic, positive
constants depending only on the dimension $q$ and other such fixed quantities in the discussion as the function $h$, the different norms
involved in the formula, etc.. Their value will be different at different
occurrences, even within the same formula. The symbol $A\sim B$ will
mean $cA \le B \le c_1A$.

The following proposition summarizes some of the important properties of the kernels defined in \eref{kerndef}.
\begin{prop}\label{kernelprop}
Let $S\ge q$ be an integer, $h:[0,\infty)\to \RR$ be  a $S$ times iterated integral of a function of bounded variation, $h(x)=1$ for $x\in [0,1/2]$, $h(x)=0$ for $x>1$, and $h$ be non--increasing. Let $\x\in\SS^q$. We have for every integer $n\ge 0$, $\Phi_n(h;\circ\cdot\x)\in \Pi_n^q$, and
\be\label{polyreprod}
\int \Phi_n(h;\x\cdot\zeta)P(\zeta)d\mu_q(\zeta)=P(\x), \qquad P\in\Pi_{n/2}.
\ee
Further,
\bea\label{kernl1bdcond}
\lefteqn{\sup_{n\ge 1,\ \zeta\in\SS^q}\int
|\Phi_n(h;\zeta\cdot\xi)|d\mu_q(\xi)=\sup_{n\ge 1}\int |\Phi_n(h;\x\cdot\xi)|d\mu_q(\xi)}\nonumber\\
&=& \frac{2\pi^{q/2}}{\Gamma(q/2)}\sup_{n\ge 1}\int_{-1}^1|\Phi_n(h;u)|(1-u^2)^{q/2-1}du<\infty,
\eea
\bea\label{kernnormest}
\lefteqn{\int |\Phi_n(h;\x\cdot\xi)|^2d\mu_q(\xi)=\frac{2\pi^{q/2}}{\Gamma(q/2)}\int_{-1}^1|\Phi_n(h;u)|^2(1-u^2)^{q/2-1}du}\nonumber\\
& \sim& n^q\sim \max_{\xi\in\SS^q} |\Phi_n(h;\x\cdot\xi)|=|\Phi_n(h;1)|.
\eea
and for every $\xi\in\SS^q$, $\xi\not=\x$,
\be\label{decaycond}
|\Phi_n(h;\x\cdot\xi)|\le cn^q\left\{\begin{array}{ll}
 (n\sqrt{1-\x\cdot\xi})^{1/2-q/2-S}, & \mbox{ if } 0\le \x\cdot\xi <1,\\
 n^{-S}, &  \mbox{ if } -1\le \x\cdot\xi <0.
\end{array}\right.
\ee
\end{prop}

Except for \eref{kernnormest}, all parts of Proposition~\ref{kernelprop} have been proved and verified repeatedly in \cite{jacobframe, locsmooth, besovpap, zfquadpap}.  We will sketch a proof of this proposition, mainly to reconcile notations.

\noindent\textsc{Proof of Proposition~\ref{kernelprop}.} 
The equation \eref{polyreprod} and the first two equations in \eref{kernl1bdcond} are clear.  The last estimate in \eref{kernl1bdcond} follows from \cite[Lemma~4.6]{jacobframe} with following choice of  the parameters there: $\alpha=\beta=q/2-1$, $h_\nu=h(\nu/n)$, where we observe that by a repeated application of the mean value theorem,
$$
\sum_{\nu=0}^\infty (\nu+1)^s|\Delta^rh(\nu/n)| \le cn^{s-r+1}, \qquad s\in\RR,\ r,n=1,2,\cdots,
$$
where $\Delta^r$ is the $r$--th order forward difference applied with respect to $\nu$.
Similarly, the estimate \eref{decaycond} follows from \cite[Lemma~4.10]{jacobframe} with same parameters as above,  $S$ in place of $K$ in \cite{jacobframe}, and $y=\x\cdot\xi$ (cf. Appendix to \cite{besovpap}).  We prove \eref{kernnormest}.  The first equation is a consequence of the rotation invariance of $\mu_q$. In view of the addition formula \eref{addformula}, 
\be\label{pf3eqn1}
\Phi_n(h;\x\cdot\xi)=\sum_{\ell=0}^n h(\ell/n)\sum_{k=1}^{d_\ell^q} Y_{\ell,k}(\x)Y_{\ell,k}(\xi).
\ee
It follows using \eref{constchristoffel}  and the facts that $h(\ell/n)=1$ for $\ell\le n/2$, $0\le h(t)\le 1$ for $t\in [0,\infty)$, that
$$
\int \Phi_n(h;\x\cdot\xi)^2d\mu_q(\xi)=\sum_{\ell=0}^n h(\ell/n)^2\sum_{k=1}^{d_\ell^q} Y_{\ell,k}(\x)^2=\sum_{\ell=0}^n h(\ell/n)^2d_\ell^q \sim n^q.
$$
Similarly, using Schwarz inequality,  \eref{addformula}, \eref{constchristoffel}, and the fact that $h(\ell/n)\ge 0$,
\begin{eqnarray*}
\lefteqn{\Phi_n(h;1)=|\Phi_n(h;\x\cdot\x)| \le \sup_{\xi\in\SS^q}
|\Phi_n(h;\x\cdot\xi)|}\\
&\le& \sum_{\ell=0}^n h(\ell/n)\{\sum_{k=1}^{d_\ell^q}Y_{\ell,k}(\x)^2\}^{1/2}\{\sum_{k=1}^{d_\ell^q}Y_{\ell,k}(\xi)^2\}^{1/2}=  \sum_{\ell=0}^n h(\ell/n)d_\ell^q=\Phi_n(h;1).
\end{eqnarray*}
Since $d_\ell^q\sim \ell^{q-1}$, $0\le h(\ell/n)\le 1$, and $h(\ell/n)=1$ for $\ell\le n/2$, the above two estimates lead to \eref{kernnormest}.
\qed

In the remainder of this paper, $h$ will denote a fixed function satisfying the conditions of Proposition~\ref{kernelprop}.

\subsection{Quadrature formulas}\label{quadexistsect}
Let $\C$ be a finite set of distinct points on $\SS^q$. A \emph{quadrature formula} based on $\C$ has the form ${\cal Q}(f)=\sum_{\xi\in\C}w_\xi f(\xi)$, where $w_\xi$, $\xi\in\C$, are real numbers. For integer $n\ge 0$, the formula is exact for degree $n$   if ${\cal Q}(P)=\int_{\SS^q}Pd\mu_q$ for all $P\in\Pi_n^q$. It is not difficult to verify that if ${\cal Q}_n(f)=\sum w_{\xi_n}f(\xi_n)$ is a sequence of quadrature formulas, with ${\cal Q}_n$ being exact with degree $n$, then ${\cal Q}_n(f)\to\int fd\mu_q$ for every continuous function $f$ on $\SS^q$ if and only if $\sum |w_{\xi_n}| \le c$, with $c$ being independent of $n$. In the sequel, we will assume tacitly that $\C$ is one of the members of a nested sequence of finite subsets of $\SS^q$, whose union is dense in $\SS^q$. All the constants may depend upon the whole sequence, but not on any individual member of this sequence. Thus, a formula ${\cal Q}$ will be called a bounded variation formula if $\sum_{\xi\in\C}|w_\xi| \le c$, with the understanding that this is an abbreviation for the concept described above with a sequence of quadrature formulas.

\begin{definition}\label{mzdef}
Let $m\ge 0$ be an integer. The set $\C$ admits an M--Z quadrature of order $m$ if there exist weights $w_\xi$ such that 
\begin{equation}\label{quadrature}
 \int_{\SS^q} P(\x)d\mu_q(\x) = \sum_{\xi\in\C}w_\xi
P(\xi), \qquad P\in \Pi_{2m}^q,
\end{equation}
and
\be\label{mzineq}
\left(\sum_{\xi\in\C}|w_\xi| |P(\xi)|^p\right)^{1/p} \le c \|P\|_p, \qquad P\in \Pi_{2m}^q, \ 1\le p< \infty.
\ee
 The weights $w_\xi$ will be called M--Z weights of order $m$. The condition \eref{mzineq} will be referred to as the M--Z condition.
\end{definition}

If $\C$ admits  an M--Z quadrature of order $m$, and $\{w_\xi\}$ are the weights involved, it is clear from using \eref{mzineq} with the polynomial identically equal to $1$ in place of $P$ that $\sum_{\xi\in\C}|w_\xi|\le c$. Further, if $\zeta\in\C$, then applying \eref{mzineq} with $p=2$ and $\Phi_m(h;\zeta\cdot\circ)$ in place of $P$, we obtain for M--Z weights of order $m$:
$$
|w_\zeta| \Phi_m(h;1)^2\le \sum_{\xi\in\C}|w_\xi| \Phi_m(h;\zeta\cdot\xi)^2\le c\int \Phi_m(h;\zeta\cdot\x)^2d\mu_q(\x).
$$
The estimate \eref{kernnormest} now implies that for all M--Z weights $\{w_\xi\}$ of order $m$,
\be\label{wtbds}
|w_\xi|\le cm^{-q}, \qquad \xi\in\C.
\ee

In \cite{mnw1}, we proved that every finite  set $\C\subset \SS^q$ admits an M--Z quadrature with an order  depending upon how dense the set $\C$ is. This density is measured in terms of the mesh norm. The mesh norm of $\C$ with respect to a subset $K\subseteq \SS^q$ 
 is defined to be
\begin{equation}
\label{meshsize}
\delta_\C(K):=\sup_{\x\in K} \mbox{\rm dist}(\x,\C).
\end{equation}

The following theorem summarizes the quadrature formula given in
\cite{mnw1}. 
\begin{theorem}\label{mzstuff}
There exists a constant $\alpha_q$ with the following
property. Let  $\C$ be a finite set of distinct
points on $\SS^q$, and $m$ be an integer with $m \le \alpha_q
(\delta_\C(\SS^q))^{-1}$. Then $\C$ admits an M--Z quadrature of order $m$, and the set $\{w_\xi\}$       of M--Z weights may be chosen to satisfy 
\be\label{no_of_nodes}
|\{\xi\ :\ w_\xi \not= 0\}| \sim m^q \sim \hbox{\rm dim}(\Pi_{2m}^q).
\ee
\end{theorem}

\bhag{Polynomial operators}\label{operatorsect}

For $t>0$, we define the summability operator  $\sigma^*_t$ by the formula
\be\label{sigmastardef}
\sigma^*_t(h;f,\x) =\int_{\SS^q} f(\zeta)\Phi_t(h;\x\cdot\zeta)d\mu_q(\zeta)=\sum_{\ell=0}^\infty h(\ell/t)\sum_{k=1}^{d_\ell^q}\hat f(\ell,k)Y_{\ell,k}(\x), \qquad f\in L^1, \ \x\in\SS^q.
\ee
(It is convenient, and customary in approximation theory, to use the notation $\sigma_t^*(h;f,\x)$ rather than $\sigma_t^*(h;f)(\x)$.) Although we defined the operator for $L^1$ to underline the fact that it is a universal operator, we will be interested only in its restriction to $C(\SS^q)$. If $f :\SS^q\to\RR$ is a continuous function, the degree of approximation of $f$ from $\Pi_x^q$ is defined by
$$
E_x(f)=\inf_{P\in\Pi_x^q}\|f-P\|_\infty.
$$
It is well known \cite{lizorkin, locsmooth} that for all integer $n\ge 1$, and $f\in C(\SS^q)$,
\be\label{goodapproxcont}
E_{n}(f)\le \|f-\sigma^*_n(h;f)\|_\infty \le cE_{n/2}(f).
\ee

Following \cite{locsmooth}, we now define a discretized version of these operators.

If $\C\subset \SS^q$ is a finite set, $\W=\{w_\xi\}_{\xi\in\C}$ and ${\bf Z}=\{z_\xi\}_{\xi\in\C}$ are sets of real numbers, we define the polynomial operator
\be\label{sigmagendef}
\sigma_t(\C, \W;h;{\bf Z},\x):=\sum_{\xi\in\C}w_\xi z_\xi \Phi_t(h;\x\cdot\xi),\qquad t\in\RR,\ \x\in\SS^q.
\ee
If $f:\SS^q\to\RR$, and $z_\xi=f(\xi)$, $\xi\in\C$, we will write $\sigma_t(\C,\W;h;f,\x)$ in place of 
$\sigma_t(\C,\W;h;{\bf Z},\x)$. In \cite{locsmooth}, we had denoted these operators by $\sigma_t(\nu;h,f)$, where $\nu$ is the measure that associates the mass $w_\xi$ with $\xi\in\C$. In this paper, we prefer to use the slightly expanded notation as in \eref{sigmagendef}. If $n\ge 1$ is an integer, $\C\subset\SS^q$ is a finite set that admits an M--Z quadrature of order $n$, $\W$ is the set of the corresponding M--Z weights.  Then it is shown in \cite[Proposition~4.1]{locsmooth} that 
\be\label{goodapprox}
E_{n}(f)\le \|f-\sigma_n(\C,\W;h;f)\|_\infty \le cE_{n/2}(f), \qquad f\in C(\SS^q). 
\ee

In this paper, we will be especially interested in the approximation of functions in the class  $\WW_r$, $r>0$, comprised of functions $f\in C(\SS^q)$ for which $E_n(f)=\O(n^{-r})$, $n\ge 1$. A complete characterization of the classes $\WW_r$ in terms of such constructive properties of its members as the number of partial derivatives and their moduli of smoothness is well known \cite{Pawelke, lizorkin}. In view of  \eref{goodapproxcont}, $f\in\WW_r$ if and only if
$$
\|f\|_{\WW_r} := \|f\|_\infty +\sup_{n\ge 1}2^{nr}\|\sigma^*_{2^n}(h;f)-\sigma^*_{2^{n-1}}(h;f)\|_\infty <\infty.
$$

 In practical applications, the data is contaminated with noise. Therefore, we wish to examine the behavior of our operators based on a data of the form $\{(\xi, f(\xi)+\epsilon_\xi)\}$, where $\epsilon_\xi$ are independent random variables with unknown probability distributions, each with mean  $0$. If the range of these random variables is not bounded, one can still assume that the probability of the variables going out of a sufficiently large interval is small. Hence, it is customary in learning theory to assume that the variables $\epsilon_\xi$ have a bounded range, so that one may use certain technical inequalities of probability theory, known as Bennett's inequalities; see the proof of Lemma~\ref{genpolystatlemma} below.   

In the statements of the theorems below, we use three parameters. The symbol $M$ denotes the number of points in the data set; we assume that the set admits an M--Z quadrature of order $m$, and the degree $n$ of the polynomial approximant $\sigma_n(\C,\W;h;f)$ is determined in terms of $m$. For theoretical considerations where one is not concerned about the actual numerical constructions of the quadrature weights, one may imagine a data set $\C$ with  $M:=|\C|\sim \delta_C(\SS^q)^{-q}$, and assume that the weights $\W$ are as guaranteed by Theorem~\ref{mzstuff}. If so, then we may take $M\sim m^q\sim \delta_\C(\SS^q)^{-q}$ in the discussion in this section. For example, the estimates \eref{soberrbd} and \eref{goodapproxnoisy} below can then be expressed in
 terms of the number of samples respectively as follows:
\be\label{goodapproxsampbd1}
 \|f-\sigma_n(\C,\W;h;f)\|_\infty\le cM^{-r/q}, \qquad \mbox{with some $n\sim M^{1/q}$}, 
\ee
and
\be\label{goodapproxsampbd2}
\prob\left(\|\sigma_n(\C,\W; h;{\bf Z})-f\|_\infty\ge c_1\frac{(\log M)^c}{M^{r/(q+2r)}}\right) \le c_2M^{-c}, \qquad \mbox{with some $n\sim (M/\log M)^{1/(2r+q)}$}.
\ee

\begin{theorem}\label{probopertheoa}
Suppose that $m\ge 1$ is an integer, $\C=\{\xi_j\}_{j=1}^M$  admits an M--Z quadrature of order $m$, and let $\W$ be the corresponding quadrature weights. Let $ r >0$, $f\in \WW_r$, $\|f\|_{\WW_r}=1$. \\
{\rm (a)} For integer $n\le m$, we have
\be\label{soberrbd}
\|f-\sigma_n(\C,\W;h;f)\|_\infty\le cn^{-r}. 
\ee
{\rm (b)}
For $j=1,\cdots, M$, let $\e_j$ be independent random variables  with mean $0$ and range $[-1,1]$,  ${\bf Z}=\{\e_j +f(\xi_j)\}$. If  $A>0$ and $n\ge 1$ is the greatest integer with $(A+q)n^{2r+q}\log n\le c_3m^q$, $n\le m$, then 
\be\label{goodapproxnoisy}
\prob\left(\|\sigma_n(\C,\W; h;{\bf Z})-f\|_\infty\ge c_1n^{-r}\right) \le c_2n^{-A}.
\ee
Here, the constants $c_1,c_2,c_3$ are independent of the distribution of the variables $\e_j$.
\end{theorem}

We now turn our attention to local approximation by our operators.    In the sequel, if $K\subseteq\SS^q$, $f:K\to\RR$, then $\|f\|_{\infty,K}:=\sup_{\x\in K}|f(\x)|$.
If $\x_0\in\SS^q$, a function $f$ is defined to be $r$--smooth at $\x_0$ if  there is a spherical cap $\SS^q_\delta(\x_0)$ such that $f\phi\in \WW_r$ for every infinitely differentiable function $\phi$ supported on $\SS^q_\delta(\x_0)$. We have proved in \cite[Theorem~3.3]{locsmooth} that $f$ is $r$--smooth at a point $\x_0$ if and only if there is a cap $\SS^q_\delta(\x_0)$ such that
$$
\|\sigma^*_{2^n}(h;f)-\sigma^*_{2^{n-1}}(h;f)\|_{\infty,\SS^q_\delta(\x_0)} =\O(2^{-nr}).
$$
Accordingly, if $K$ is a spherical cap, we may define the class $\WW_r(K)$ to consist of $f\in C(\SS^q)$, for which
$$
\|f\|_{\WW_r(K)}:=\|f\|_\infty + \sup_{n\ge 1}2^{nr} \|\sigma^*_{2^n}(h;f)-\sigma^*_{2^{n-1}}(h;f)\|_{\infty,K} <\infty.
$$

\begin{theorem}\label{localapproxtheo}
Suppose that $m\ge 1$ is an integer, $\C=\{\xi_j\}_{j=1}^M$  admits an M--Z quadrature of order $m$, and let $\W$ be the corresponding quadrature weights. Let $0<r\le S-q$, $K'\subset K$ be concentric spherical caps,  $f\in C(\SS^q)$, $\|f\|_{W_r(K)}=1$.\\
{\rm (a)} For integer $n$, $1\le n \le m$,
\be\label{localapproxnonoise}
\|f-\sigma_n(\C,\W;h;f)\|_{\infty, K'} \le cn^{-r}.
\ee 
{\rm (b)} For $j=1,\cdots, M$, let $\e_j$  be independent random variables  with mean $0$ and range contained in $[-1,1]$, and ${\bf Z}=\{\e_j +f(\xi_j)\}$. If  $A>0$ and $n\ge 1$ is the greatest integer with $(A+q)n^{2r+q}\log n\le c_3m^q$, $n\le m$, then 
\be\label{goodlocapproxnoisy}
\prob\left(\|\sigma_n(\C,\W; h;{\bf Z})-f\|_{\infty,K'}\ge c_1n^{-r}\right) \le c_2n^{-A}.
\ee
Here, the constants $c_1,c_2,c_3$ are independent of the distribution of the variables $\e_j$.
\end{theorem}. 

\bhag{Construction of quadrature formulas}\label{quadconsect}
In this section, we describe two algorithms to obtain bounded variation quadrature formulas associated with a given finite set of points $\C\subset \SS^q$. Both of these constructions can be described in a very general setting. Since this also simplifies the notations and ideas considerably by avoiding the use of real and imaginary parts of a doubly indexed polynomial $Y_{\ell,k}$, we will describe the algorithms in this generality.

Let $\Omega$ be a nonempty set, $\mu$ be a probability measure on $\Omega$, $\C\subset \Omega$, $y_1,y_2,\cdots$ be a complete orthonormal basis for $L^2(\Omega,\mu)$,  where $y_1\equiv 1$, and $V_k$ denote the span of $y_1,\cdots,y_k$.  Let $\nu$ be another measure on $\Omega$, and $\ip{\circ}{\circ}$ denote the inner product of $L^2(\Omega,\nu)$. For integer $N\ge 1$, the Gram matrix $G_N$ is an $N\times N$ matrix, defined by $(G_N)_{\ell,k}=\ip{y_\ell}{y_k}=(G_N)_{k,\ell}$, $1\le k,\ell\le N$. We wish to find a weight function $W$ on $\Omega$ such that $\int_\Omega Pd\mu=\int_\Omega PWd\nu$ for all $P\in V_N$ for an integer $N$ for which $G_N$ is positive definite. 

For the applications to the case of quadrature formulas for the sphere,
$\Omega=\SS^q$, $\mu=\mu_q$,  and $y_k$'s are the orthogonal spherical
harmonics, arranged in a sequence, so that $y_1\equiv 1$, and all
polynomials of lower degree are listed before those of a higher degree.
To include all polynomials in $\Pi_n^q$, we need $N=d_n^{q+1}$. 
There are many possibilities to define the measure $\nu$. The simplest is the measure $\nu^{MC}$ that associates the mass $1/|\C|$ with each point of $\C$. A more sophisticated way to define the measure $\nu$ is the following. We obtain a partition of $\SS^q$ into a dyadic triangulation such that each triangle contains at least one point of $\C$. We choose only one point in each triangle, and hence, assume that each triangle contains exactly one point of $\C$. We define the measure $\nu^{TR}$ to be the measure that associates with each $\xi\in\C$ the area of the triangle containing $\xi$. 

One of the simplest ideas to compute the quadrature weights is the following. Let $N$ be an integer for which $G_N$ is positive definite. If $P=\sum_{j=1}^N a_jy_j$, then $\int_\Omega Pd\mu=a_1$. Also, the vector ${\bf a}=(a_1,\ldots,a_N)^T$ satisfies the matrix equation
$$
G_N{\bf a}=(\ip{P}{y_1},\ldots,\ip{P}{y_N})^T,
$$
so that
\be\label{abstractquad1}
\int_\Omega Pd\mu=\sum_{k=1}^N(G_N)^{-1}_{1,k}\ip{P}{y_k}= \ip{P}{\sum_{k}(G_N)^{-1}_{1,k}y_k}.
\ee
In the setting of the sphere,  this gives the following quadrature formula:
\be\label{sphquad1}
\int_{\SS^q}Pd\mu_q=\sum_{\xi\in\C} P(\xi)\left\{\nu(\{\xi\}) \sum_{k=1}^N(G_N)^{-1}_{1,k}y_k(\xi)\right\}=: \sum_{\xi\in\C}w_\xi^{LSQ}P(\xi).
\ee
We formulate this as \\

\leftline{\textsc{Algorithm LSQ}}
\textbf{Input:} The matrix $Y=(y_k(\xi))$,
$k=1,\cdots,N$ (optional), and the vector ${\bf v}=(\nu(\{\xi\})$. 
\begin{enumerate}
\item\label{lsqsolveforb} Solve $Y\texttt{diag}({\bf v})Y^T{\bf b}=(1,0,\ldots,0)^T$. 
\item\label{returnwts} Return $w_\xi^{LSQ} = \nu(\{\xi\})\sum_{k=1}^N b_ky_k(\xi)$.
\end{enumerate}

We observe that $G_N=Y\texttt{diag}({\bf v})Y^T$. It is clear that the
matrix $G_N$ is always positive semi--definite; the assumption that it
is positive definite is equivalent to the assumption that no element of
$V_N$  vanishes identically on $\C$. If $\C$ and
$\{\nu(\{\xi\})\}$ satisfy the M--Z inequalities,
Theorem~\ref{lsqwttheo} below shows that $G_N$ is well conditioned. Assuming that the matrix $Y$ is input, the time to compute $G_N$ is $\O(N^2|\C|)$ and the space requirement is  $\O(N^2)$. (In the case of the sphere $\SS^q$, we need $N=d_n^{q+1}=\O(n^q)$ to compute formulas exact for degree $n$.) The vector ${\bf b}$ in
Step~\ref{lsqsolveforb} can be found using such iterative methods  as
the conjugate residual method.   We refer to  \cite{bernd_fischer} for a more detailed analysis of this method.
Using this approach, the matrix $Y$ and $G_N$ need not  be stored or 
precomputed, but the product of the matrix $G_N$ with an arbitrary residual 
vector ${\bf r}$ needs to be computed. This observation results in a substantial saving in the time and memory complexity of the algorithm when the results are desired only within a given accuracy. For the unit sphere $\SS^2$, 
when $N=(n+1)^2$, the product $G_N {\bf r} =Y\texttt{diag}({{\bf v}})Y^T {\bf r}$ can be computed within an accuracy $\e$ using a recent algorithm 
of Keiner \cite{keiner} using $\O(n^2(\log n)^2 + \log(1/\epsilon)|\C|)$
operations, where $\epsilon$ is the accuracy of the method.  

One way to interprete this algorithm is the following. Let $f :\SS^q\to\RR$, 
and $P$ be  the solution to the least square problem
$$
P=\arg\min\{\ip{f-Q}{f-Q}\ :\ Q\in V_N\}.
$$
If ${\bf f}$ is the vector $(\ip{f}{y_j})$, then $P=\sum_j (G_N^{-1}{\bf f})_j y_j$. The quadrature formula with weights $w_\xi^{LSQ}$ thus offers $\int_{\SS^q}Pd\mu_q$ as the approximation to $\int_{\SS^q}fd\mu_q$. 
The weights $w_\xi^{LSQ}$ also satisfy a least square property among all the possible quadrature formulas, as shown in  Lemma~\ref{lsqwtlemma}(a). We summarize some of the properties of the weights $w_\xi^{LSQ}$ in the following theorem. 

\begin{theorem}\label{lsqwttheo}
Let $n\ge 1$ be an integer, $N=d_n^{q+1}$, $\C$ be a finite set of points on $\SS^q$ and $\nu$ be a measure supported on $\C$. Let $v_\xi:=\nu(\{\xi\})$, $\xi\in \C$, and
\be\label{vmz}
c_1\|P\|_p\le \left\{\sum_{\xi\in\C} v_\xi |P(\xi)|^p\right\}^{1/p}\le c_2 \|P\|_p, \qquad P\in \Pi_n^q,\ 1\le p \le \infty,
\ee
{\rm (a)} For the Gram matrix $G_N$, the lowest eigenvalue is $\ge c_1^2$, and the largest eigenvalue is $\le c_2^2$, where $c_1,c_2$ are the constants in \eref{vmz} with $p=2$. In particular, $G_N$  is positive definite. Moreover, $\sum_{\xi\in\C}|w_\xi^{LSQ}|\le c$.\\
{\rm (b)} If 
\be\label{strongmz}
\left|\int P^2d\nu -\int P^2d\mu_q\right| \le \frac{c}{n^q}\int P^2d\mu_q, \qquad P\in\Pi_{n}^q,
\ee
then $|w_\xi^{LSQ}|\le cv_\xi$, $\xi\in\C$. In particular, the weights $\{w_\xi^{LSQ}\}$ satisfy the M--Z condition.\\
{\rm (c)} Let $M\ge 1$ be an integer, $\C$ be a independent random sample of $M$ points chosen from the distribution $\mu_q$, and $A, \eta>0$. Let $v_\xi =1/M$, $\xi\in\C$. There exists a constant $c=c(A)$ such that if $n\ge 2$ is an integer with $M\ge cn^q\log n/\eta^2$, then 
\be\label{pottsineq}
\prob\left(\left|\int P^2d\nu -\int P^2d\mu_q\right| \ge \eta\int P^2d\mu_q, \qquad P\in\Pi_{n}^q\right)\le c_1n^{-A}.
\ee
In particular, if $M\ge cn^{3q}\log n$, then the condition \eref{strongmz} is satisfied with probability exceeding $1-c_1n^{-A}$.
\end{theorem}

One disadvantage of the algorithm LSQ is that one needs to know the
value of $N$ in advance. We now describe an idea which has the potential to avoid this problem. In the case when $\Omega$ is a subset of a Euclidean space, 
and the $y_j$'s are polynomials, with $y_1$ denoting the constant polynomial, one can
construct a system $\{t_k\}$ of orthonormalized polynomials with respect to $\nu$
using recurrence relations. 
Recurrence relations for orthogonal polynomials in
several variables have been discussed in detail by Dunkl and Xu \cite[Chapter~3]{xubook}. 
In contrast to the viewpoint in \cite{xubook}, we may depend upon a specific enumeration, 
but require the recurrence relation to have a specific form described in Theorem~\ref{recurtheo} below. 
This form allows us to generalize the ideas in Gautschi's book \cite[Chapter~2]{gautschi} in our context.

To describe our ideas in general, let $\Omega\subset \RR^{q+1}$, $u_1,u_2,\cdots$ be the lexicographic enumeration of the monomials in $q+1$ variables, so that $u_1$ is the monomial identically equal to $1$, the restrictions of $u_k$'s to $\Omega$ are linearly independent, and $V_k=\mbox{span }\{u_1,\cdots,u_k\}$. It is not difficult to see that for every integer $k\ge 1$, there is a minimal index $p(k)$ such that there exists a monomial $\tilde f_k$ of degree $1$ with
\be\label{monomialrec}
\tilde f_ku_{p(k)}=u_{k+1}, \qquad k=1,2,\cdots.
\ee
We now let, for each $k=1,2,\cdots$, $\{y_1,\cdots,y_k\}$   be a basis for $V_k$ orthonormal with respect to $\mu$, $N\ge 1$ be an integer for which the Gram matrix $G_N$ is positive definite, and for each $k=1,\cdots,N$, $\{t_1,\cdots,t_k\}$ be a basis for $V_k$ orthonormal with respect to $\nu$. 
Clearly,  any polynomial $P\in V_N$ can be written in the form
$$
P(\x)=\int P(\zeta)\sum_k t_k(\x)t_k(\zeta)d\nu(\zeta),
$$
and consequently, one gets the ``quadrature formula''
\be\label{absrecquad}
\int P(\x)d\mu(\x)=\int P(\zeta) \left\{\sum_k \left(\int t_k(\x)d\mu(\x)\right) t_k(\zeta)\right\}d\nu(\zeta).
\ee
In this discussion only, let 
$t_k=:\sum_{j}c_{k,j}y_j$, and the matrix $(c_{k,j})$ be denoted by $C$. The condition that $t_1,\cdots,t_N$ is an orthonormal
system with respect to $\nu$ is equivalent to the condition that $CG_NC^T=I$, where $I$ is the $N\times N$ identity matrix. 
Hence, $G_N^{-1}=C^TC$. Moreover, $\int t_kd\mu=c_{k,1}$ for $k=1,\cdots,N$, and hence, we conclude that
$$
\sum_k \left(\int t_k(\x)d\mu(\x)\right) t_k=\sum_j\sum_k c_{k,1}c_{k,j}y_j =\sum_j (G_N)^{-1}_{1,j}y_j.
$$
Thus, the quadrature weights in \eref{absrecquad} are the same as those in \eref{abstractquad1}.

First, we summarize the various recurrence relations in Theorem~\ref{recurtheo} below, although we will not use all of them.
We will denote the (total) degree of $u_k$  by $D_k$, and observe that $D_k$ is also the degree of $y_k$ and $t_k$, $D_j\le j$,  and $D_{p(k)}=D_{k+1}-1$.

\begin{theorem}\label{recurtheo}
There exist real numbers $s_{k,j}$, $\tilde r_{k,j}$,  $A_k\ge 0$, and a linear polynomial $f_k$, such that
\be\label{ytrecur}
f_ky_{p(k)}=y_{k+1} -\sum_{D_{k+1}-2\le D_j\le D_{k}\atop j\le k} \tilde r_{k,j}y_j, \  f_kt_{p(k)}= A_kt_{k+1}-\sum_{D_{k+1}-2\le D_j\le D_{k}\atop j\le k}s_{k,j}t_j.
\ee
More generally, if $P$ is any linear polynomial, there exist real numbers $r_{k,j}(P)$ such that
\be\label{yrecurgen}
 P y_k = \sum_{D_k-1\le D_j\le D_k+1} r_{k,j}(P)y_j, 
\ee
We have $t_k=\sum_{j}c_{k,j}y_j$, where
\be\label{connectrecur}
A_kc_{k+1,\ell}=\left\{\sum_{D_\ell-1\le D_m\le D_\ell+1} r_{\ell,m}(f_k)c_{p(k),m} + \sum_{D_{k+1}-2\le D_j\le D_{k}}s_{k,j}c_{j,\ell}\right\}.
\ee
\end{theorem}

In the context of the sphere $\SS^q$, we will compute $t_k$'s using \eref{ytrecur}, and compute $\int t_kd\mu_q$ using a 
known quadrature formula. The resulting algorithm, Algorithm REC, in the context of the sphere is summarized below. 
This algorithm is similar to the Stieltjes method in Gautschi's book \cite[Section~2.2]{gautschi}. Even though it is 
feasible to carry out the algorithm for as large an $N$ as the data allows, and to find this value of $N$ during run time, 
it is still desirable from the point of view of numerical stability to limit the largest $N$ from the outset. Accordingly, 
in describing the following algorithm, we stipulate that the quadrature formula is to be computed to be exact only for 
polynomials in $V_N$ for the largest possible $N\le L$ for some integer 
$L\ge 1$. We assume further that we know another quadrature formula (for example, the Driscoll--Healy formula \cite{driscoll_healy}) exact for polynomials in $V_L$:
\be\label{seedquad}
\sum_{\zeta\in \C^*} \lambda_\zeta P(\zeta) =\int Pd\mu_q, \qquad P\in V_L.
\ee

\noindent
\leftline{\textsc{Algorithm REC}}
\textbf{Input:} An integer $L$, the sequence $p(k)$, $k=1,\cdots, L$, sets $\C$, $\C^*$, weights $(\lambda_\zeta)_{\zeta\in\C^*}$ so that \eref{seedquad} holds,   the values $\{y_j(\xi)\}_{\xi\in\C}$, $\{y_j(\zeta)\}_{\zeta\in \C^*}$ for $j=1,2,3,4$, and the values $f_k(\xi)$, $f_k(\zeta)$, $k=1,\cdots,L$.
\begin{enumerate}
\item Using Gram--Schmidt procedure, initialize $t_1, t_2,t_3,t_4$, both for points in $\C$ and in $\C^*$, and initialize $N=4$. 
\item For $k=1,\cdots,4$, let $\gamma_k=\sum_{\zeta\in\C^*} \lambda_\zeta t_k(\zeta)$. 
\item For each $\xi\in\C$, initialize $w_\xi = \sum_{k=1}^4\gamma_kt_k(\xi)$.
\item For $k= 4, 5,\cdots$ (so that the degrees are at least $0$ for all polynomials entering in the recursions) and while $N\le L$, repeat steps \ref{initrecur}--\ref{lastrecur} below.
\item\label{initrecur} For $j$ with $D_{k+1}-2\le D_j\le D_{k}$, set
$$s_{k,j}=\ip{f_kt_{p(k)}}{t_j}.$$
\item\label{trecur} Define $T_{k+1}$ by 
$$
T_{k+1}=f_kt_{p(k)}-\sum_{D_{k+1}-2\le D_j\le D_{k}}s_{k,j}t_j.
$$
both for points in $\C$ and points in $\C^*$. If $I_{k+1}=\ip{T_{k+1}}{T_{k+1}}=0$, then stop, and set $N=k$. Otherwise, define $t_{k+1}=T_{k+1}/I_{k+1}^{1/2}$.
\item Set $\gamma_{k+1}=\sum_{\zeta\in\C^*} \lambda_\zeta t_{k+1}(\zeta)$.
\item\label{lastrecur} For each $\xi \in\C$, $w_\xi=w_\xi+\gamma_{k+1}t_{k+1}(\xi)$, $k=k+1$, $N=N+1$.
\end{enumerate}

In the case of the sphere $\SS^q$, we take   $L=d_{\tilde n}^{q+1}$ for some integer ${\tilde n}\ge 1$.  The number of $j$'s with $D_{k+1}-2\le D_j\le
D_k$, $1\le j,k\le L$ is $\O({\tilde n}^{q-1})$. In this discussion only, let $M=|\C|+|\C^*|$.
Consequently, Steps~\ref{initrecur} and \ref{trecur} require $\O(M{\tilde n}^{q-1})$ operations. Since the remaining two steps in the
loop take $\O(M)$ operations, the loop starting at Step~4 require $\O(M{\tilde n}^{2q-1})$ operations. Finally, we observe that in
implementing the above algorithm, one need not keep  the whole matrix $t_k(\xi)$; only the rows corresponding to three degrees
are required in any step. In particular, the memory requirement of this
algorithm is $\O(M{\tilde n}^{q-1})$.

\bhag{Numerical experiments}\label{numsect}
The objective of this section is to demonstrate and supplement the theoretical results presented in 
Sections~\ref{operatorsect} and \ref{quadconsect}. 

Our first set of experiments illustrates the algorithms LSQ and REC. The experiments were conducted over a long period of time, many of them long before we started to write the paper. Hence, the normalizations for the spherical polynomials $Y_{\ell,k}$ are somewhat different in Tables~\ref{lsqtable}  and \ref{tableREC} from the rest of the paper. This is reflected in the sum of the absolute values of the weights, but has no effect on the various results other than scaling.

First, we report on the algorithm LSQ. Each of the experiments in this case was repeated $30$ times with data sets 
chosen randomly from the distribution $\mu_2$ on $\SS^2$. To test our algorithms, we computed the computed Gram matrix 
$G^{COM}$ given by
$$
G^{COM}_{\ell,m}=\sum_{\xi\in\C}w_\xi^{LSQ} y_\ell(\xi)y_m(\xi), \quad \ell,m
               < \lfloor n/2 \rfloor.
$$
The average maximum matrix norm of the difference between $G^{COM}$ and the identity matrix of the same size indicates 
the error of the quadrature formulas. The results are shown in Table~\ref{lsqtable}.    Based on these results we conjecture that in order to obtain stable quadrature formulas (i.e., with small condition number for the original Gram matrix $G_N$) exact for degree $n\ge 1$, one has to use at most $4d_n^{q+1}$ uniformly distributed points. In contrast, the theoretical guarantee in Theorem~\ref{lsqwttheo}(c) requires $\O(n^{3q}\log n)$ points.
\begin{table}[h]
$$
\begin{array}{|c|c|c|c|c|c|c|c|c|c|}\hline
M & n & Error & \sum |w_\xi| & \min w_\xi & \max w_\xi & \mbox{pos}  
                                       & \kappa(G_N) & \lambda_{\min} &  \lambda_{\max} \\
\hline
8192& 16& 2.41*10^{-15} & 3.5449 &  2.29*10^{-4}& 7.88*10^{-4} & 8192 &  2.43   & 0.607  & 1.4730 \\
    & 44& 4.32*10^{-15} & 3.5714 & -5.06*10^{-4}& 0.0029       & 8039 &  37.52  & 0.078  & 2.8047 \\
    & 64& 6.15*10^{-15} & 5.5575 & -0.00664     & 0.0073       & 6068 &  1695.1 & 0.003  & 3.9315 \\
    & 84& 9.73*10^{-12} & 82.152 & -0.16274     & 0.1551       & 4431 &  3.52*10^6& 2.6*10^{-6} & 5.4851\\ 
\hline
16384&44 &4.43*10^{-15}& 3.5449 & -9.50*10^{-6} & 8.82*10^{-4} &16382   & 9.19   & 0.24036 & 2.1590\\
     &64 &5.25*10^{-15}& 3.5787 & -3.75*10^{-4} & 0.0015 &16014 & 51.9   & 0.05964 & 2.9150\\ 
     &84 &7.10*10^{-15}& 4.4757 & -0.0024    & 0.0032    &13361 & 944.86 & 0.00612 & 3.8457\\
     &100&1.94*10^{-15}& 9.1325 & -0.0077    & 0.0072    &10625 & 11896.1& 4.8*10^{-4}& 4.6008\\
\hline
32768& 44 &6.02*10^{-15} & 3.5449&  3.11*10^{-5}& 2.90*10^{-4}& 32768 & 4.270 & 0.4157 & 1.7652 \\
     & 64 &7.09*10^{-15} & 3.5450& -1.79*10^{-5}& 5.23*10^{-4}& 32761 & 7.977 & 0.8208 & 5.2519 \\
     & 84 &7.71*10^{-15} & 3.5574& -1.43*10^{-4}& 7.92*10^{-4}& 32410 & 42.97 & 0.0716 & 2.7967 \\
     &100 &7.62*10^{-15} & 3.6777& -4.28*10^{-4}& 9.96*10^{-4}& 30819 & 145.6 & 0.0250 & 3.2967 \\
\hline
\end{array}
$$
\caption{The statistics for the experiments with the algorithm LSQ. $M=|\C|$,
$n-2$ is the degree of spherical polynomials for which exact quadrature formulas
were computed, $N=n^2$, pos stands for the number of positive weights, $\kappa(G_N)$, $\lambda_{\min}$, $\lambda_{\max}$ are the condition number, the maximum eigenvalue and the minimum eigenvalue of the matrix $G_N$ respectively.}
\label{lsqtable}
\end{table}

As can be seen from the table, for a fixed degree $n$, the condition number 
$\kappa(G_N)$ decreases as the number of points increases. For $n>140$ and various sets of randomly generated points 
on the sphere, we do not obtain good numerical results. 
This might be due to  a defect in the built in numerical procedures used by Matlab in computing the spherical harmonics  of high degree at values close to $-1$ or $1$. The situation was much better for the dyadic points; i.e., the centers of the dyadic triangles.

For dyadic points on the sphere, the best result we obtained so far is $n=178$ with $131,072$ points.  
As a further verification of this quadrature, we considered the following data. The data are constructed 
using coefficients $\{a_{\ell,k}\}$ for spherical polynomials up to degree $90$, taken from model MF4 
used for modelling the lithospheric field. The model is computed by geophysicists at 
GeoForschungsZentrum Potsdam (Germany) based on CHAMP satellite data. We use those coefficients to 
construct the samples of a function $f=\sum_{\ell,k}a_{\ell,k}Y_{\ell,k}$ at the centers of $8*4^7$ 
dyadic triangles. We then use our pre-computed quadrature based at these centers which can integrate 
spherical polynomials up to degree $178$ to compute the Fourier coefficients $\widehat{a}_{\ell,k}$. 
The maximum difference between the vector $\{\widehat{a}_{\ell,k}\}$ and the vector $\{a_{\ell,k}\}$ 
was found to be $6.66*10^{-15}$.

Next, we considered the algorithm REC. In the context of spherical polynomials, the recurrence relations have to be chosen very carefully using 
the special function properties of the spherical harmonics $Y_{\ell,k}$, in order to get stable 
results \cite{PottsSteidlTasche}. In the present situation, the polynomials $t_k$ have no special structure. 
Therefore, it turns out that the algorithm REC  is not very stable for high degrees. However, when we took the centers 
of $8192$ dyadic triangles as the quadrature nodes, and used the measure $\nu^{TR}$ as the starting measure, then we are able to obtain satisfactory quadrature formulas for degree $32$. We note an interesting feature here that all the weights obtained by this algorithm are positive. 
These results are summarized in Table~\ref{tableREC} below.
\begin{table}[h]
$$
\begin{array}{|c|c|c|c|c|}\hline
 n  &  \mbox{Error} &  \min(w_\xi) &  \max(w_\xi) &  \sum w_\xi \\
\hline
16 &  4.196643*10^{-14} & 5.181468*10^{-4} & 2.538441*10^{-3} & 12.56637 \\
     22 &  5.302425*10^{-13} & 5.175583*10^{-4} & 2.543318*10^{-3} & 12.56637 \\
     32 &  9.240386*10^{-11} & 5.154855*10^{-4} & 2.544376*10^{-3} & 12.56637 \\
     42 &  4.434868*10^{-8}  & 5.086157*10^{-4} & 2.544141*10^{-3} & 12.56637 \\
     44 &  2.320896*10^{-5}  & 5.094771*10^{-4} & 2.562948*10^{-3} & 12.56638 \\ 
\hline
\end{array}
$$
\caption{Quadrature constructed using REC on $8192$ dyadic points}
\label{tableREC}
\end{table}

Our second set of experiments demonstrates the local approximation properties of the operators $\sigma_n(\C,\W;h)$ for a smooth function $h$. For this purpose, we consider the following benchmark functions, considered by various authors \cite{SlWo:04, sloanrbf, psdiff, ganmhas1}, listed in \eref{benchfns} below. Using the notation $\x=(x_1,x_2,x_3)$,
the functions are defined by
\bea\label{benchfns}
g_1(\x) &=& (x_1-0.9)^{3/4}_{+} + (x_3-0.9)^{3/4}_{+},\nonumber\\
g_2(\x) &=& [0.01 - (x_1^2+x_2^2+(x_3-1)^2)]_{+} + \exp(x_1+x_2+x_3),\nonumber \\
g_3(\x) &=& 1/(101-100x_3),\nonumber\\
g_4(\x) &=& 1/(|x_1|+|x_2|+|x_3|),\nonumber\\
g_5(\x) &=&\left\{\begin{array}{ll} \cos^2\left(\frac{3\pi}{2}{\rm dist}(\x,(-1/2,-1/2, 1/\sqrt{2}))\right),
& {\rm {if~~dist}}(\x,(-1/2,-1/2, 1/\sqrt{2})) < 1/3,  \\
0, &  {\rm {if~~dist}}(\x,(-1/2,-1/2, 1/\sqrt{2})) \geq 1/3. 
\end{array}\right.
\eea
In order to define the function $h$, we recall first that the $B$ spline $B_m$ of order $m$ is defined recursively \cite[p.~131]{deboor} by
\be\label{bsplinedef}
B_m(x):=\left\{\begin{array}{ll}
1, &\mbox{ if $m=1$, $0<x\le 1$,}\\
0, &\mbox{ if $m=1$, $x\in\RR\setminus (0,1]$,}\\
\disp\frac{x}{m-1}B_{m-1}(x)+\frac{m-x}{m-1}B_{m-1}(x-1), &\mbox{ if $m>1$, $x\in\RR$.}
\end{array}\right.
\ee
The function $B_m$ is an $m-1$ times interated integral of a function of bounded variation. We will choose $h$ to be
\be\label{hfuncdef}
h_m(x)=\sum_{k=-m}^{m} B_{m}(2mx-k),
\ee
for different values of $m$, in order to illustrate the effect of the smoothness of $h_m$ on the quality of local approximation. If $m\ge 3$, the function $h_m$ satisfies the conditions in Proposition~\ref{kernelprop} with $S=m-1$. We note that the discretized Fourier projection operator $\sigma_{63}(\C,\W;h_1)$ has been called the hyperinterpolation operator \cite{sloan}.

One example of the localization properties of our operators is given in the following table, where we show the error in approximation of $g_1$ on the whole sphere and on the cap $K=\SS^2_{0.4510}((-1/\sqrt{2},0,-1/\sqrt{2})$. The operators were constructed using Driscoll--Healy quadrature formula \cite{driscoll_healy} based on $4(n+1)^2$ points, exact for integrating polynomials of degree $2n$. The maximum error on the whole sphere, given in Columns~2 and 3, is estimated by the error at $10000$ randomly chosen points; that on the cap, given in Columns~4 and 5, is estimated by the error at $1000$ randomly chosen points on the cap. It is clear that even though the maximum error on the whole sphere is slightly better for the (discretized) Fourier projection than for our summability operator, the singularities of $g_1$ continue to dominate the error in the Fourier projection on a cap away from these singularities; the performance of our summability operator is far superior.
\begin{table}[ht]
$$
\begin{array}{|c|c|c|c|c|}
\hline
\mbox{n}& S2errh1 & S2errh5 &Kerrh1 & Kerrh5\\
\hline
63 &     0.0097&                 0.0112&  3.4351*10^{-4}&  6.5926*10^{-7}\\
\hline
  127 &    0.0044&                 0.0055&     8.0596*10^{-5} &  6.5240*10^{-8}\\
\hline
  255 &    0.0033&                 0.0038& 1.4170*10^{-5} &  1.1816*10^{-8}\\
\hline
\end{array}
$$
\caption{$S2errh1=\max_{\x\in\SS^2}|g_1(\x)-\sigma_n(\C,\W;h_1,g_1,\x)|$, $S2errh5=\max_{\x\in\SS^2}|g_1(\x)-\sigma_n(\C,\W;h_5,g_1,\x)|$, $Kerrh1=\max_{\x\in K}|g_1(\x)-\sigma_n(\C,\W;h_1,g_1,\x)|$, $Kerrh5=\max_{\x\in K}|g_1(\x)-\sigma_n(\C,\W;h_5,g_1,\x)|$, $(\C,\W)$ are given by the Driscoll--Healy formulas. }
\label{maxerrtab}
\end{table} 

Theorem~\ref{localapproxtheo} points out another way to demonstrate the superior localization of our summability operator without an a priori knowledge of the locations of the singularities. Since each of the test functions is infinitely differentiable on large caps of different sizes, Theorem~\ref{localapproxtheo} suggests that the more localized the method, the greater is the probability that the approximation error would be smaller than a given number.  To demonstrate also how our ability to construct quadrature formulas based on scattered data helps us to analyse the approximation properties of our summability operators, we took for  the set $\C$ a randomly generated sample of $65536$ points. For these points, the  weights $\W$ computed by the algorithm LSQ yield a quadrature formula  exact for integrating spherical polynomials of degree $126$. 
We compare three approximation methods, the least square approximation from $\Pi_{63}^2$, the approximation given by the operator $\sigma_{63}(\C,\W;h_1)$, and the approximation given by $\sigma_{63}(\C,\W;h_5)$.  For each function, we computed the absolute value of the difference between the approximate value computed by each of the three methods and the true value of the function at $20,000$ randomly chosen points on the sphere.  The percentage of points where the value of this difference is less than $10^{-x}$ is reported in Table~\ref{optesttable} below, for $x=2:10$. It is very obvious that $\sigma_{63}(\C,\W;h_5)$ gives a far superior performance than the other methods, due to its localization properties. 

\begin{table}[ht]
$$
\begin{array}{|c|c|c|c|c|c|c|c|c|c|c|}
\hline
&x\rightarrow & 10&9&8&7&6&5&4&3&2\\
\hline
&S1&    0  &   0&     0.005&  0.02&   0.42&   4.44&  39.43&  94.79&  100 \\
\cline{2-11}
g_1&LS&    0&     0&     0&      0.04&   0.56&   5.32&  46.38&  95.45&  100\\
\cline{2-11}
&S5&    0.02&  0.19&  1.87&  16.89&  59.36&  68.34&  79.01&  93.09&   99.97 \\
\hline
&S1&    0&     0.01&   0.09&  0.74&   7.94&  84.99&  99.19&  99.99&  100 \\
\cline{2-11}
g_2&LS&    0&     0.01&   0.15&  1.29&  13.29&  86.09&  99.28&  100&    100 \\
\cline{2-11}
&S5&    0.39&  3.34&  41.95& 90.78&  94.52&  97.19&  99.18&  99.97&  100 \\
\hline
&S1&    0&     0.01&   0.11&  1.26&  12.02&  91.87&  99.86&  100&    100 \\
\cline{2-11}
g_3&LS&    0&     0.01&   0.10&  1.49&  16.26&  93.12&  99.87&  100&    100 \\
\cline{2-11}
&S5&    0.51&  5.43&  51.08& 82.22&  91.90&  95.79&  98.49&  99.87&  100 \\
\hline
&S1&    0&     0&      0&     0.01&   0.18&   1.91&   18.28&  83.81&   99.97 \\
\cline{2-11}
g_4&LS&    0&     0&      0.01&  0.02&   0.25&   2.16&   21.24&  86.43&   99.98 \\
\cline{2-11}
&S5&    0.01&  0.01&   0.04&  0.36&   3.47&  17.48&   40.98&  80.06&   99.88 \\
\hline
&S1&    0&     0&    0.01&    0.09&   1.12&   11.84&  88.94&  99.45&  100\\
\cline{2-11}
g_5&LS&    0&     0.01& 0.01&    0.15&   1.42&   15.23&  90.47&  99.75&  100\\
\cline{2-11}
&S5&    0.08&  0.64& 5.73&   66.82&  83.54 &  88.74&  92.95&  96.78&   97.64\\
\hline
\end{array}
$$
\caption{Percentages of error less than $10^{-x}$ for different functions, LS= Least square, S1= error with $\sigma_{63}(\C,\W;h_1)$, S5= error with $\sigma_{63}(\C,\W;h_5)$. For example, for the function $g_3$, S5 was less than $10^{-7}$ for 82.22\% of the 20000 randomly selected points, while S1 (respectively, LS) was less than $10^{-7}$ for 1.26\% (respectively, 1.49\%) points.}
\label{optesttable}
\end{table}

\medskip\noindent
 Next, we illustrate the stability of our operators under noise. Since our operators are linear operators, we assume for this part of the study that the target function $f$ is the zero function
contaminated either by uniform random noise in the range
$[-\epsilon,\epsilon]$, or a normally distributed random variable with
mean $0$ and standard deviation $\e$. We let $\C$ be a set of $65536$ random points and computed corresponding weights $\W$ that integrate exactly polynomial up to degree $126$. 
These were used in calculating $\sigma_{63}(\C,\W;h_1)$ and
$\sigma_{63}(\C,\W;h_5)$ at each point of a test data set consisting of $20000$ random samples from the distribution $\mu_2$.   For each value
of $\epsilon = 0.1; 0.01; 0.001; 0.0001$, the experiment is repeated $50$ times
and the errors are the averaged over the number of repetitions.  The percentage of points at which the absolute computed  value is less
than $10^{-x}$ is reported in Table \ref{noisetable} in the case when $\epsilon=0.01$. The results for the other values of $\epsilon$ were consistent with the linearity of the operator.
We observe that in each case, both $\sigma_{63}(\C,\W;h_1)$ and $\sigma_{63}(\C,\W;h_5)$ yield  better results than the least squared approximation, while $\sigma_{63}(\C,\W;h_5)$ is slightly superior to $\sigma_{63}(\C,\W;h_1)$.

\begin{table}
$$
\begin{array}{|c|c|c|c|c||c|c|c|c|}
\hline
   x\rightarrow &   5& 4& 3& 2&3.0& 2.75&2.5&    2.25 \\  \hline 
S1&       0.05&    0.635& 9.93&  97.45& 0&     7.97&   92.75&   100.00   \\     
LS&          0&        0&    0&    100&     0&     0.04&   30.93&    99.07   \\    
S5&       0.085&   1.015& 10.03& 97.49&  0.24&    51.97&   99.87&   100.00\\   
\hline
\end{array}
$$
\caption{Percentages of error less than $10^{-x}$ for $\epsilon=0.01$,
LS= Least square, S1= error with $\sigma_{63}(\C,\W;h_1)$, S5= error
with $\sigma_{63}(\C,\W;h_5)$. The random noise in the left half comes from the uniform 
distribution in $[-\epsilon,\epsilon]$, that in the right half from the normal distribution with mean $0$, standard deviation $\e$.}
\label{noisetable}
\end{table}

\begin{figure}[t]
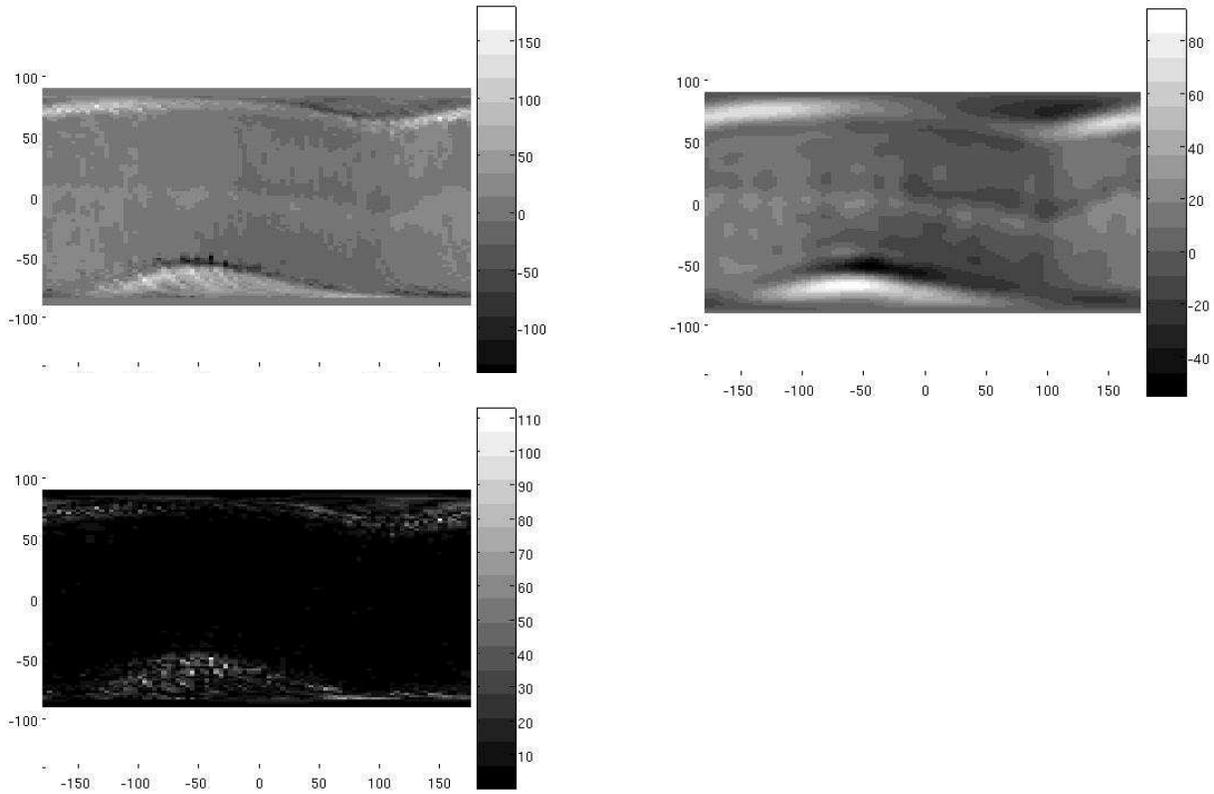

\begin{minipage}[t]{0.45\textwidth}
 \includegraphics[width=\textwidth]{newfYOlsen.epsi}
\end{minipage}
\begin{minipage}[t]{0.45\textwidth}
 \includegraphics[width=\textwidth]{newL23q7_Y.epsi}
\end{minipage}
\begin{minipage}[t]{0.45\textwidth}
  \includegraphics[width=\textwidth]{newL23q7_err.epsi}
\end{minipage}
\caption{From left to right: The original data, its reconstruction using $\sigma_{22}(\C,\W;h_7)$, and the error in the approximation, $|\sigma_{22}(\C,\W;h_7)- y|$.}
\label{magsatfig}
\end{figure}

Finally, we used our operator $\sigma_{22}(\C,\W;h_7)$ with the MAGSAT data. Our purpose here is only to test how our methods work on a ``real life'' data. This data, supplied
to us kindly by Dr. Thorsten Maier, measures the magnetic field of the earth in nT as a vector 
field. It was derived from vectorial MAGSAT morning data that has been processed by Nils Olsen
of the Danish Space Research Institute. The measurements are averaged on a longitude-latitude grid
with $\Delta \phi = 4^{o}$ and $\Delta \theta = 2^{o}$ in geomagnetic coordinates. The radial
variations of the MAGSAT satellite have been neglected in the dataset and, therefore, prior to the
averaging process, the GSFC(12/83) reference potential model has been subtracted. The data
results from one month of measurements, centered at March 21, 1980. We extract the East West component of the vectorial data as a scalar valued
function on the sphere.  Totally, there are $8190$ 
data sites. 
A quadrature of degree $44$ was computed based on those sites. Figure~\ref{magsatfig} shows
the original data, its reconstruction using $\sigma_{22}(\C,\W;h_7)$, and the error in the
approximation, $|\sigma_{22}(\C,\W;h_7)- y|$, as a map in the longitude-latitude plane. 
As can be seen from the figures, the reconstruction preserves the key features of the original data.

\bhag{Proofs}\label{proofsect}
In the interest of organization, we will prove the various new results in the
paper in the following order. We will prove  Theorem~\ref{recurtheo} first, since its proof does not require any preparation. We will then use Proposition~\ref{kernelprop} to prove Theorems~\ref{probopertheoa}(a) and \ref{localapproxtheo}(a).  Next, we will prove Lemma~\ref{lsqwtlemma} and use it to prove parts (a) and (b) of Theorem~\ref{lsqwttheo}. The  remaining results in this paper involve probabilities. We prove   Lemma~\ref{genpolystatlemma} next, estimating the probability that the supremum norm of a sum of random spherical polynomials exceeds a given number.  This  lemma will be  used immediately to prove Theorem~\ref{lsqwttheo}(c). Finally, we will prove  Theorems~\ref{probopertheoa}(b), \ref{localapproxtheo}(b).

\noindent
\textsc{Proof of Theorem~\ref{recurtheo}.} 
It is convenient to prove \eref{yrecurgen} first. Since $Py_k$ is a polynomial of degree $D_k+1$, there exist real numbers $r_{k,j}(P)$ such that
$$
Py_k=\sum_{D_j\le D_k+1} r_{k,j}(P)y_j.
$$
Since the system $\{y_k\}$ is orthonormal with respect to $\mu$, 
$$
r_{k,j}(P)=\int_{\Omega} Py_ky_jd\mu.
$$
If $D_j<D_k-1$, then the degree of $Py_j$ is less than $D_k$. Because of the lexicographic ordering where lower degree polynomials appear before the higher degree ones, this implies that $Py_j\in V_{k-1}$. Since $y_k$ is orthogonal to $V_{k-1}$, it follows that $r_{j,k}(P)=0$ if $D_j<D_k-1$. This completes the proof of \eref{yrecurgen}.

We observe that $y_{p(k)}\in \mbox{span }\{u_1,\cdots,u_{p(k)}\}$. So, there exists a constant $\alpha$ such that $\alpha y_{p(k)}- u_{p(k)}\in V_{p(k)-1}$. Thus, $\alpha\tilde f_ky_{p(k)}-\tilde f_ku_{p(k)}=\alpha\tilde f_ky_{p(k)}-u_{k+1}$ is linear combination of terms of the form $\tilde f_k u_j$, $1\le j\le p(k)-1$.  Since $p(k)$ is the minimal index for which there exists a monomial $\tilde f_k$ with $\tilde f_k u_{p(k)}\in V_{k+1}$, each of the terms $\tilde f_k u_j$, $1\le j\le p(k)-1$ is in $V_k$. It follows that $\alpha\tilde f_ky_{p(k)}-u_{k+1}\in V_k$. Again, there exists a constant $\alpha'$ such that $\alpha'u_{k+1}-y_{k+1}\in V_k$. Therefore, writing $f_k=\alpha\alpha' \tilde f_k$, we conclude that $f_ky_{p(k)}-y_{k+1}=\alpha'(\alpha\tilde f_ky_{p(k)}-u_{k+1})+\alpha'u_{k+1}-y_{k+1}\in V_k$; i.e., 
$\disp
f_ky_{p(k)}=y_{k+1}-\sum_{D_j\le D_k}\tilde{r}_{k,j}y_j. 
$
The first equation in \eref{ytrecur} is now proved in view of \eref{yrecurgen}, applied with $p(k)$ in place of $k$, and the fact that $D_{p(k)}=D_{k+1}-1$. We note that $f_k$ is a constant  multiple of the monomial $\tilde f_k$. The  second equation in \eref{ytrecur} is proved in the same way. 

Using the second equation in \eref{ytrecur} and \eref{yrecurgen}, we obtain from the definition of $c_{k,j}$'s that
\begin{eqnarray*}
A_kc_{k+1,\ell}&=&A_k\int_\Omega t_{k+1}y_\ell d\mu\\
&=&  \int_\Omega f_kt_{p(k)}y_\ell d\mu+\sum_{D_{k+1}-2\le D_j\le D_{k}}s_{k,j}\int_\Omega t_jy_\ell d\mu\\
&=& \sum_{D_\ell-1\le D_m\le D_\ell+1} r_{\ell,m}(f_k)\int_\Omega t_{p(k)}y_m d\mu +\sum_{D_{k+1}-2\le D_j\le D_{k}}s_{k,j}c_{j,\ell}\\
&=& \sum_{D_\ell-1\le D_m\le D_\ell+1} r_{\ell,m}(f_k)c_{p(k),m} + \sum_{D_{k+1}-2\le D_j\le D_{k}}s_{k,j}c_{j,\ell}.
\end{eqnarray*}
This proves \eref{connectrecur}.
\qed

Next, we use Proposition~\ref{kernelprop} to prove Theorem~\ref{probopertheoa}(a) and Theorem~\ref{localapproxtheo}(a).

\noindent
\textsc{Proof of Theorem~\ref{probopertheoa}}(a). To prove part (a), we assume without loss of generality that $n\ge 8$, and let $\ell\ge 1$ be the largest integer with $2^{\ell+2}\le n$. In view of \eref{goodapprox},
\begin{eqnarray*}
\|f-\sigma^*_n(\C,\W;h;f)\|_\infty &\le& cE_{n/2}(f)\le cE_{2^{\ell+1}}(f)\le c\|f-\sigma_{2^{\ell+1}}^*(h;f)\|_\infty\\
&\le& c\sum_{k=\ell+1}^\infty \|\sigma_{2^{k+1}}^*(h;f)-\sigma_{2^{k}}^*(h;f)\|_\infty\le c2^{-r\ell}\le cn^{-r}.
\end{eqnarray*}
This proves part (a). \qed

\noindent
\textsc{Proof of Theorem~\ref{localapproxtheo}}(a). 
 Let $K''$ be a spherical cap, concentric with $K$, $K'$, and having radius equal to the average of the radii of $K$, $K'$. Let $\psi$ be fixed, $C^\infty$ function that is equal to $1$ on $K''$ and equal to $0$ outside of $K$. Without loss of generality, we may assume that $n\ge 8$, and let $\ell\ge 1$ be the largest integer such that $2^{\ell+2}\le n$. The direct theorem of approximation theory (cf. \cite{Pawelke}) implies that there exists $P\in \Pi_{2^\ell}^q$ such that
$$
\|\psi-P\|_\infty \le c2^{-\ell S}.
$$
Therefore, using the definition of $\|f\|_{\WW_r(K)}$, we conclude that
\begin{eqnarray*}
E_{2^{\ell+1}}(f\psi)&\le& \|f\psi-P\sigma^*_{2^\ell}(h;f)\|_{\infty}\le  \|(f-\sigma^*_{2^\ell}(h;f))\psi\|_{\infty} +\|(\psi-P)\sigma^*_{2^\ell}(h;f)\|_{\infty}\\
&\le& c\{\|f-\sigma^*_{2^\ell}(h;f)\|_{\infty,K} + 2^{-nS}\|f\|_{\infty}\}\\
&\le& \disp c\left\{\sum_{k=\ell+1}^\infty \|\sigma_{2^{k+1}}^*(h;f)-\sigma_{2^{k}}^*(h;f)\|_{\infty,K} +2^{-nS}\|f\|_{\infty}\right\}\le c2^{-r\ell}.
\end{eqnarray*}
In view of \eref{goodapprox}, 
\bea\label{pf4eqn1}
\|f-\sigma_n(\C,\W;h;f\psi)\|_{\infty, K'}&=&\|f\psi-\sigma_n(\C,\W;h;f\psi)\|_{\infty, K'}\le \|f\psi -\sigma_n(\C,\W;h;f\psi)\|_\infty\nonumber\\
&\le&  cE_{n/2}(f\psi)\le  E_{2^{\ell+1}}(f\psi)\le c2^{-r\ell}\le cn^{-r}.
\eea
Since $1-\psi(\zeta)=0$ for $\zeta\in K''$, we may use \eref{decaycond} to deduce that for $\x\in K'$,
\begin{eqnarray*}
\lefteqn{|\sigma_n(\C,\W;h;(1-\psi)f,\x)| =\left|\sum_{\xi\in\C\setminus K''}w_\xi f(\xi)(1-\psi(\xi)) \Phi_n(h;\x\cdot\xi)\right|}\\
&\le& \frac{c(K,K',K'')}{n^{S-q}}\|(1-\psi)f\|_\infty \sum_{\xi\in\C}|w_\xi| \le \frac{c(K,K',K'')}{n^{S-q}}.
\end{eqnarray*}
Together with \eref{pf4eqn1} and the fact that $r\le S-q$, this implies \eref{localapproxnonoise}. \qed

Next, we prove Lemma~\ref{lsqwtlemma}, describing certain extremal properties for the weights $w^{LSQ}_\xi$. These  will be used in the proof of parts (a) and (b) of Theorem~\ref{lsqwttheo}. 

\begin{lemma}\label{lsqwtlemma}
 Let $n\ge 1$ be an integer, $N=d_n^{q+1}$, $\C$ be a finite set of points on $\SS^q$ and $\nu$ be a measure supported on $\C$. Let $v_\xi:=\nu(\{\xi\})$, $\xi\in \C$.\\
{\rm (a)} If the Gram matrix is positive definite, then the weights $w_\xi^{LSQ}$ are solutions of the extremal problem to minimize $\sum_{\xi\in\C} {w_\xi^2}/{v_\xi}$ subject to the conditions that $\sum_{\xi\in\C}w_\xi  y_\ell(\xi)=\delta_{1,\ell}$.\\
{\rm (b)} If \eref{vmz} holds, there exist real numbers $W_\xi$, $\xi\in\C$, such that $|W_\xi|\le v_\xi$ for $\xi\in\C$ and $\sum_{\xi\in\C}W_\xi P(\xi)= \int Pd\mu_q$ for all $P\in \Pi_n^q$.
\end{lemma}
\begin{Proof} 
In this proof, we will write $G$ in place of $G_N$.
The Lagrange multiplier method to solve the minimization problem sets up parameters $\lambda_\ell$ and minimizes
$$
\sum_{\xi\in\C} {w_\xi^2}/{v_\xi}-2\sum_\ell \lambda_\ell\left(\sum_\xi w_\xi  y_\ell(\xi)-\delta_{1,\ell}\right).
$$
Setting the gradient (with respect to $w_\xi$) equal to $0$, we get $w_\xi=v_\xi\sum_\ell \lambda_\ell y_\ell(\xi)$. Writing, in this proof only, $Q=\sum_\ell \lambda_\ell y_\ell$, we see that $w_\xi=v_\xi Q(\xi)$. Substituting back in the linear constraints, this reduces to $\sum_{\xi\in\C} v_\xi Q(\xi)y_\ell(\xi)=\delta_{1,\ell}$. These conditions determine $Q$ uniquely; indeed, $Q=\sum_j G^{-1}_{1,j}y_j$. This proves part (a).

The part (b) is proved essentially in \cite{mnw1}, but since it is not stated in this manner, we sketch a proof again. During this proof, different constants will retain their values. Let $M=|\C|$, $\RR^M$ be equipped with the norm $\||{\bf r}\||=\sum_{\xi\in\C} v_\xi |r_\xi|$. In this proof only, let ${\cal S}$ be the operator defined on $\Pi_n^q$ by ${\cal S}(P)=(P(\xi))_{\xi\in\C}\in\RR^M$, and ${\mathbb V}$ be the range of ${\cal S}$. The estimate 
\be\label{pf2eqn1}
\int |P|d\mu_q\le c_1\sum_{\xi\in\C}v_\xi |P(\xi)|
\ee
 implies that the operator ${\cal S} : \Pi_n^q\to {\mathbb V}$ is invertible. We may now define a linear functional on ${\mathbb V}$  by
$$
x^*({\bf r}) =\int {\cal S}^{-1}({\bf r})d\mu_q, \qquad {\bf r}\in {\mathbb V}.
$$
It is clear from \eref{pf2eqn1} that the norm of $x^*$ is bounded above by $c_1$. The Hahn--Banach theorem yields a norm preserving extension of this functional to the whole space $\RR^M$. Identifying this functional with the vector $(W_\xi)_{\xi\in\C}$, the extension property implies that
$\sum_{\xi\in\C}W_\xi P(\xi)= \int Pd\mu_q$ for all $P\in \Pi_n^q$, while the norm preservation property implies that $|W_\xi|\le c_1v_\xi$ for $\xi\in\C$. 
\end{Proof}

\noindent
\textsc{Proof of Theorem~\ref{lsqwttheo}} (a), (b). 
Let $N=d_n^{q+1}$, ${\bf r}\in\RR^N$, and $P=\sum_{\ell}r_\ell y_\ell$. In this proof only, we write $G$ in place of $G_N$. Then
$$
{\bf r}^TG{\bf r} =\sum_{\ell,m} r_\ell \left\{\sum_{\xi\in\C}v_\xi y_\ell(\xi)y_m(\xi)\right\}r_m =\sum_{\xi\in\C} v_\xi P(\xi)^2,
$$
and ${\bf r}^T{\bf r}=\|P\|_2^2$. 
Therefore, \eref{vmz} with $p=2$ implies that $c_1^2{\bf r}^T{\bf r}\le {\bf r}^TG{\bf r} \le  c_2^2{\bf r}^T{\bf r}$ for all ${\bf r}\in\RR^N$. The statements about the eigenvalues of $G$ are an immediate consequence of  the Raleigh--Ritz theorem \cite[Theorem~4.2.2]{horn}. Using Lemma~\ref{lsqwtlemma}(b), we  obtain weights $W_\xi$ such that $\sum_{\xi\in\C}W_\xi y_\ell(\xi)=\delta_{1,\ell}$, $\ell=1,\cdots,N$, and $|W_\xi|\le cv_\xi$, $\xi\in\C$. During the remainder of this proof, we write $w_\xi=w_\xi^{LSQ}$. In view of Lemma~\ref{lsqwtlemma}(a), we have
$$
\sum_{\xi\in\C}|w_\xi|\le \left\{\sum_\xi v_\xi\right\}^{1/2}\left\{\sum_{\xi\in\C}\frac{w_\xi^2}{v_\xi}\right\}^{1/2} \le \left\{\sum_\xi v_\xi\right\}^{1/2}\left\{\sum_{\xi\in\C}\frac{W_\xi^2}{v_\xi}\right\}^{1/2}\le c\sum_{\xi\in\C}v_\xi \le c_1.
$$
This completes the proof of part (a).

In order to prove part (b), we adopt the following notation during this proof only. Let $I$ denote the $N\times N$ identity matrix. For any $N\times N$ matrix $H$, let $\|H\|$ denote $\sup\|H{\bf r}\|$, $\|{\bf r}\|=1$, ${\bf r}\in\RR^N$. We note that $\|H\|$ is the largest singular value of $H$. If $H$ is a symmetric, positive definite matrix, then it is also the largest eigenvalue of $H$, and moreover, $|{\bf r}_1^TH{\bf r}_2|\le \|H\|\|{\bf r}_1\|\|{\bf r}_2\|$, ${\bf r}_1, {\bf r}_2\in \RR^N$. Using \eref{strongmz}, it is easy to conclude using the Raleigh--Ritz theorem that $\|G-I\|\le cn^{-q}$, $\|G^{-1}\|\le c$, and hence, 
$$
\|G^{-1}-I\|=\|G^{-1}(I-G)\|\le c\|G^{-1}\|\|G-I\| \le cn^{-q}.
$$
Let $\y(\x)$ denote the vector $(y_1(\x), \cdots, y_N(\x))^T$ for $\x\in\SS^q$. In view of the addition formula, $\|\y(\x)\|^2$ is independent of $\x$, and hence, 
$$
\|\y(\x)\|^2=\int_{\SS^q} \sum_{j=1}^N y_j(\x)^2d\mu_q(\x) =d_n^{q+1}\le cn^q, \qquad \x\in\SS^q.
$$
 Consequently,  we have
\begin{eqnarray*}
\lefteqn{\frac{|w_\xi^{LSQ}|}{v_\xi}= \left|\int \y(\x)^T G^{-1}\y(\xi)d\mu_q(\x)\right|}\\
& \le& \left|\int \y(\x)^T(G^{-1}-I)\y(\xi)d\mu_q(\x)\right| + \left|\int \y(\x)^T\y(\xi)d\mu_q(\x)\right|\\
& \le& \|G^{-1}-I\|\int \|\y(\x)\|\|\y(\xi)\|d\mu_q + |(1,0,\cdots,0)^T\y(\xi)| \le cn^{-q}n^q +c\le c.
\end{eqnarray*}
This completes the proof of part (b). \qed

The proof of the remaining new results in the paper are based on the following lemma, that gives a recipe for estimating the probabilities involving polynomial valued random variables. 

\begin{lemma}\label{genpolystatlemma} 
Let $n, M\ge 1$ be  integers, $\{\omega_j\}_{j=1}^M$ be independent random variables, and for $j=1,\cdots,M$, $Z_j=Z(\omega_j,\circ)\in\Pi_n^q$ have mean equal to $0$ according to $\omega_j$. Let  $B, R>0$,  $\disp \max_{1\le j\le M,\ \x\in\SS^q}|Z_j(\x)|\le Rn^q$, and the sum of the variances of
$Z_j$ be bounded by $Bn^q$ uniformly on $\SS^q$. If $A>0$ and  $12R^2(A+q)n^q\log n\le B$  then
\be\label{genpolystatest}
\prob \left(\left\|\sum_{j=1}^M Z_j\right\|_\infty\ge \sqrt{12B(A+q)n^q\log n}\right)\le c_1n^{-A}.
\ee
Here, the  positive constant $c_1$ is independent of $M$ and the distributions of $\omega_j$. 
\end{lemma}
\begin{Proof}   The proof depends upon  Bennett's inequality \cite[p.~192]{pollardbook}. In this proof  only, we adopt a slightly different meaning for the symbols $L$, $V$, $\eta$. Let $L, V, \eta$ be  positive numbers, and $X_j$, $j=1,\cdots,M$, be independent random variables. According to Bennett's inequality, if the mean of each $X_j$ is $0$, the range of each $X_j$ is a subset of $[-L,L]$, and $V$ exceeds the sum of the variances of $X_j$,  then for $\eta>0$,
\be\label{bennetineq}
\prob \left(|\sum_{j=1}^M X_j|\ge\eta\right) \le 2\exp\left(-\frac{V}{L^2}g(L\eta/V)\right),
\ee
where, in this proof only, $g(t) := (1+t)\log(1+t)-t$. We observe that $g(t)=\int_0^t\int_0^u (1+w)^{-1}dwdu$. Therefore, if $0\le t\le 1/2$, then for $0\le w\le u\le t$, $(1+w)^{-1}\ge 2/3$, and hence, $g(t)\ge t^2/3$. Consequently, if $L\eta\le V/2$, then
\be\label{pf1eqn1}
 \prob \left(|\sum_{j=1}^M X_j|\ge\eta\right) \le 2\exp(-\eta^2/(3V)).
\ee
Now, let $\x\in\SS^q$. We apply \eref{pf1eqn1} with $Z_j(\x)$ in place of $X_j$, $Rn^q$ in place of $L$, $Bn^q$ in place of $V$, $\eta =\sqrt{3B(A+q)n^q\log n}$. Our condition on $n$ ensures that $L\eta/V\le 1/2$ with these choices. Therefore,
\be\label{pf1eqn2}
  \prob \left(|\sum_{j=1}^M Z_j(\x)|\ge\sqrt{3B(A+q)n^q\log n}\right) \le 2n^{-A-q}.
\ee
Next, in the proof only, let $P^*=\sum_{j=1}^M Z_j$, $\x^*\in\SS^q$ be chosen so that $|P^*(\x^*)|=\|P^*\|_\infty$, $\C\subset \SS^q$ be chosen so that $|\C|\sim cn^q$ and $\delta_\C(\SS^q)\le 1/(2n)$. Then we may find $\xi^*\in\C$ such that $\mbox{\rm dist}(\x^*,\xi^*)\le 1/(2n)$. Since $P^*\in\Pi_n^q$, its restriction to the great circle through $\x^*$ and $\xi^*$ is a trigonometric polynomial of order at most $n$. In view of the Bernstein inequality for these polynomials \cite[Chapter~4, (1.1)]{devlorbk}, 
$$
|P^*(\xi^*)-P^*(\x^*)|\le n\|P^*\|_\infty \mbox{\rm dist}(\xi^*,\x^*)\le (1/2)|P^*(\x^*)|.
$$
We deduce that
$$
\left\|\sum_{j=1}^M Z_j\right\|_\infty \le 2\max_{\x\in\C}\left|\sum_{j=1}^M Z_j(\x)\right|.
$$
Therefore, the event $\left\|\sum_{j=1}^M Z_j\right\|_\infty\ge 2\sqrt{3B(A+q)n^q\log n}$ is a subset of the union of the $|\C|$ events $|\sum_{j=1}^M Z_j(\x)|\ge\sqrt{3B(A+q)n^q\log n}$, $\x\in\C$. Hence, the estimate \eref{pf1eqn2} implies \eref{genpolystatest} with $c_1=2c$.
\end{Proof}

We are now in a position to prove Theorem~\ref{lsqwttheo}(c).\\

\noindent \textsc{Proof of Theorem~\ref{lsqwttheo}}(c).
 Let $\x\in\SS^q$. In this proof only, let $Z_\xi=\Phi_{4n}(h;\x\cdot\xi)-\int \Phi_{4n}(h;\x\cdot\zeta)d\mu_q(\zeta)$. Then the mean of each $Z_\xi$ is $0$, and its variance can be estimated by
$$
\int Z_\xi^2d\mu_q(\xi)\le \int (\Phi_{4n}(h;\x\cdot\xi))^2d\mu_q(\xi)\le cn^q.
$$
Finally, $|Z_\xi|\le cn^q$ for each $\xi$. Hence, we may use Lemma~\ref{genpolystatlemma} with $cM$ in place of $B$, and $c$ in place of $R$, to conclude that
$$
\prob \left(\sup_{\x\in \SS^q}\left|\frac{1}{M}\sum_{\xi\in \C}\Phi_{4n}(h;\x\cdot\xi)-\int \Phi_{4n}(h;\x\cdot\zeta)d\mu_q(\zeta)\right|\ge c_2\sqrt{\frac{n^q\log n}{M}}\right) \le cn^{-A},
$$
and with $M\ge cn^{q}\log n/\eta^2$, 
$$
\prob \left(\sup_{\x\in \SS^q}\left|\frac{1}{M}\sum_{\xi\in \C}\Phi_{4n}(h;\x\cdot\xi)-\int \Phi_{4n}(h;\x\cdot\zeta)d\mu_q(\zeta)\right|\ge \eta\right) \le cn^{-A}.
$$
Since any $P\in\Pi_{2n}^q$ can be written in the form
$$
P(\zeta)=\int P(\x)\Phi_{4n}(h;\x\cdot\zeta)d\mu_q(\x),
$$
we see that with probability exceeding $1-cn^{-A}$,
\begin{eqnarray*}
\lefteqn{\left|\frac{1}{M}\sum_{\xi\in\C}P(\xi) -\int P(\zeta)d\mu_q(\zeta)\right|}\\
&=& \left| \frac{1}{M}\sum_{\xi\in\C}\int P(\x)\Phi_{4n}(h;\xi\cdot\x)d\mu_q(\x) \right.\\
&&\qquad -\left.\int \int P(\x)\Phi_{4n}(h;\x\cdot\zeta)d\mu_q(\x)d\mu_q(\zeta)\right|\\
&\le&\int |P(\x)|\left|\frac{1}{M}\sum_{\xi\in\C}\Phi_{4n}(h;\xi\cdot\x)-\int \Phi_{4n}(h;\x\cdot\zeta)d\mu_q(\zeta)\right|d\mu_q(\x)\\
&\le& \eta\int |P(\x)|d\mu_q(\x).
\end{eqnarray*}
For $P\in\Pi_n^q$, we may now apply this estimate with $P^2\in\Pi_{2n}^q$.
\qed

Another immediate consequence of Lemma~\ref{genpolystatlemma}  is the following lemma, describing the the probabilistic behavior of the operator $\sigma_n(\C,\W;h)$.
\begin{lemma}\label{noiselemma}
Suppose that $m\ge 1$ is an integer, $\C=\{\xi_j\}_{j=1}^M$  admits an M--Z quadrature of order $m$, and let $\W$ be the corresponding quadrature weights.  Let $R, V>0$, and for $j=1,\cdots, M$, $\e_j$  be independent random variables  with mean $0$, variance not exceeding $V$, and range $[-R,R]$.   Let $g\in C(\SS^q)$, $\|g\|_\infty\le 1$, and ${\bf E}=\{\e_jg(\xi_j)\}_{\xi_j\in\C'}$. Then  for integer $n\ge 1$ with $(R^2/V)(A+q)n^q\log n\le c_3m^q$, 
\be\label{opstochastic}
\prob \left(\left\|\sigma_n(\C,\W;h;{\bf E})\right\|_\infty
\ge
c_1\sqrt{\frac{V(A+q)n^q\log n}{m^q}}\right) \le c_2n^{-A}.
\ee
Here, the positive constants $c_1$, $c_2$, $c_3$  depend only on   $q$ but not on $M$ and the distributions of $\e_j$.
\end{lemma}

\begin{Proof} In this proof only, if $\xi=\xi_j\in\C$, we will write $\e_\xi$ for $\e_j$ and $w_\xi$ for the weight in $\W$ corresponding to $\xi$. We use Lemma~\ref{genpolystatlemma} with  $E_\xi=m^qw_\xi \e_\xi g(\xi)\Phi_n(h;\xi\cdot\circ)$, $\xi\in\C$. We note that the random variable $\omega_j$ in Lemma~\ref{genpolystatlemma} is $\e_\xi$ in this case. It is clear that the mean of each $E_\xi$ is $0$. Since \eref{wtbds} implies that $|w_\xi|\le cm^{-q}$, \eref{kernnormest} shows that $\|E_\xi\|_\infty \le cRn^q$. Moreover, for any $\x\in \SS^q$, the variance of $E_\xi(\x)$ does not exceed $Vm^{2q}w_\xi^2\Phi_n(h;\xi\cdot\x)^2$. In view of the fact that $w_\xi$ are M--Z quadrature weights,  \eref{wtbds} and \eref{kernnormest} imply that
\begin{eqnarray*}
\lefteqn{\sum_{\xi\in\C}m^{2q}w_\xi^2\Phi_n(h;\xi\cdot\x)^2 \le cm^q\sum_{\xi\in\C}|w_\xi|\Phi_n^2(h;\xi\cdot\x)}\\
&\le& cm^q\int_{\SS^q}\Phi_n^2(h;\zeta\cdot\x)d\mu_q(\zeta)\le cm^qn^q.
\end{eqnarray*}
Thus, we may choose $B$ in Lemma~\ref{genpolystatlemma} to be $cVm^q$. The estimate \eref{opstochastic} now follows as a simple consequence of \eref{genpolystatest}. 
\end{Proof}

We are now in a position to prove the probabilistic assertions of Theorems~\ref{probopertheoa} and \ref{localapproxtheo}.

\noindent
\textsc{Proof of Theorem~\ref{probopertheoa}}(b).
To prove part (b), we use Lemma~\ref{noiselemma} with $g\equiv 1$. Since the range of $\e_j$'s is contained in $[-1,1]$, we may take $R=V=1$, and obtain from \eref{opstochastic} that
$$
 \prob \left(\left\|\sigma_n(\C,\W;h;{\bf E})\right\|_\infty
\ge
c_4\sqrt{\frac{(A+q)n^q\log n}{m^q}}\right) \le c_2n^{-A}.
$$
The choice of $n$ with an appropriate $c_3$, ensures that $(A+q)n^q\log n/m^q \le n^{-2r}$. Therefore,
$$
 \prob \left(\left\|\sigma_n(\C,\W;h;{\bf E})\right\|_\infty
\ge
c_4n^{-r}\right) \le c_2n^{-A}.
$$
The estimate \eref{goodapproxnoisy} is now clear in view of \eref{goodapprox} and the linearity of the operators $\sigma_n(\C,\W;h)$.
\qed

\noindent
\textsc{Proof ofTheorem~\ref{localapproxtheo}}(b). 
We apply Lemma~\ref{noiselemma} again with $g\equiv 1$. As before, we may choose $R=V=1$. The choice of $n$ with an appropriate $c_3$, ensures that $(A+q)n^q\log n/m^q \le n^{-2r}$. 
Therefore, \eref{opstochastic} with these choices implies that
$$
\prob \left(\left\|\sigma_n(\C,\W;h;{\bf E})\right\|_\infty
\ge
c_4n^{-r}\right) \le c_2n^{-A}.
$$
Together with \eref{localapproxnonoise} and the linearity of the operators $\sigma_n(\C,\W;h)$, this leads to \eref{goodlocapproxnoisy}.
\qed

\bhag{Conclusion} 
We have described a construction of linear operators yielding spherical polynomial approximations based on scattered data on a Euclidean sphere. While the operators can be defined for arbitrary continuous functions on the sphere, without any a priori knowledge about the location and nature of its singularities, they are auto--adaptive in the sense that the approximation properties of these globally defined polynomials adapt themselves on the different parts of the sphere according to the smoothness of the target function on these parts. While the theoretical properties of these operators and their localization were studied in \cite{locsmooth}, a bottleneck in their numerical construction was the construction of quadrature formulas based on scattered data, exact for integrating moderately high degree spherical polynomials. So far, it was possible only to compute quadrature formulas exact at most for degree 18 polynomials. We show that a simple--minded construction involving a Gram matrix is surprisingly well conditioned, and yields the necessary quadrature rules, up to degree $178$. Using these newly constructed quadrature formulas, we are able to demonstrate that our constructions yield superior approximation properties to those of more traditional techniques of least squares and Fourier projection, in the sense that the presence of singularities in some parts of the sphere  affects the degree of approximation by our operators on other parts far less than in the case of these other traditional techniques. We give probabilistic estimates on the local and global degrees of approximation by our operators in the presence of noise, and demonstrate its use in the modeling of a ``real life'' data set. We also describe a theoretical algorithm for construction of data dependent multivariate orthogonal polynomials and their use in the construction of quadrature formulas, analogous to the univariate algorithms in the book \cite{gautschi} of Gautschi.



\end{document}